%% file: main.tex
\documentclass{article}
\usepackage[utf8]{inputenc}
\usepackage[left=2cm, right=2cm,bottom=3cm]{geometry}
\usepackage[round]{natbib}
\usepackage{hyperref}
\usepackage{amsmath}
\usepackage{amsfonts}
\usepackage{amssymb}
\usepackage[]{algorithm2e}
\RestyleAlgo{boxruled}
\usepackage{multirow}
\LinesNumbered
\usepackage[short]{optidef}
\usepackage{physics}
\usepackage[caption=false]{subfig}

\usepackage{pgfplots}
\usetikzlibrary{calc}

\DeclarePairedDelimiter\floor{\lfloor}{\rfloor}

\usepackage{xcolor}
\newcommand{\colr}[1]{{\color{black} {#1}}}

\title{Timeslot allocation for waiting list control}
\author{Y.M. van der Vlugt, J.T. van Essen, R.F.M. Vromans, M. Carlier }

\begin{document}

\maketitle

\input{Abstract}
\input{introduction}
\input{literature}

\input{model}
\input{solution_methods}

\input{data_analysis}

\section{Computational results}\label{sec:results}

\input{results}

\input{results_smk}

\input{conclusion}

\bibliographystyle{plainnat}
\bibliography{library}

\appendix

\input{abbreviations}

\end{document}

%% file: Abstract.tex
\begin{abstract}
\colr{As pressure on the healthcare system increases, patients that require elective surgery experience longer access times to pre- and post-operative appointments and surgery. Hospitals can control their waiting lists by allocating timeslots to the different types of appointments they discern. To allow appointments to be planned timely, they need to decide this allocation several weeks in advance.} However, the precise consequences of \colr{the timeslot allocation are} uncertain, as not all patients follow the same treatment pathway. Furthermore, as these planning decisions are made far in advance, they are based on an uncertain prediction of future waiting lists. \colr{We aim to develop methods that support hospitals to improve their timeslot allocation} to reduce access times for patients and ensure that all available capacity in the outpatient department and operating room is used.

The problem is modelled as a Markov decision process (MDP). As the state space is very large, an exact solution \colr{cannot be determined. Therefore, we use} least-squares policy iteration to find an approximate solution, formulate an (integer) linear program which is used to solve a deterministic variant of the MDP, and investigate \colr{several} decision rules. The solution methods are tested on a case study at the Sint Maartenskliniek, a hospital focusing on orthopaedic care in Nijmegen, the Netherlands. Based on a simulation study, we find that all methods improve on the static allocation method currently used by the hospital, with the (integer) linear program leading to the best results. 
However, the performance deteriorates with the number of weeks the hospital plans ahead. To counter this, we propose a method in which a percentage of timeslots is statically allocated far in advance, and the remaining timeslots are allocated when enough information is available to effectively deal with variability. For the case study of one surgeon at the SMK, we find that statically allocating 60\% of the timeslots and dynamically allocating the remainder 6 weeks in advance provides the best results in terms of meeting access time targets and efficient resource utilization.
\end{abstract}

%% file: introduction.tex
\section{Introduction}\label{sec:introduction}
The COVID-19 pandemic has shown in an unprecedented way how critical healthcare capacity can be to the well-being of society. 
Even in non-pandemic times, the Dutch healthcare system has been facing increasing pressure due to an ageing population and shortage of hospital staff \citep{NZa2020}. This can lead to longer access times, lower quality of care, higher costs, and overworked staff. One way to avoid such capacity shortages is to invest more in human and material resources. Another way is to make more efficient use of the resources already available by optimizing healthcare resource allocation.

\colr{Resource allocation plays an important role in scheduling outpatient appointments at the outpatient clinic for surgical specialists. Before surgery, elective} patients have one or more outpatient appointments. During these appointments, the surgeon determines the nature of the complaint, whether the patient \colr{should and can} undergo surgery, and which type of surgery is required. \colr{After surgery, a patient might have one or more follow-up appointments, and typically, a discharge appointment.} Hospitals discern between different types of appointments \colr{to control access of different types of patients, to deal with differences in urgencies, and to accomodatedifferent appointment durations, such as longer timeslots for patients that the surgeon sees for the first time.} 
It is desirable to schedule these appointments in such a way that patient access times remain stable and as many patients as possible receive treatment. 

In this paper, we explore how \colr{allocating a number of timeslots to each appointment type} can be used to control waiting lists with the objective to reduce access times of late patients and make efficient use of resources. 
\colr{We refer to} determining the number of each type of appointment to be scheduled as the \textit{timeslot allocation problem}. Our aim is to model this problem, and find and implement appropriate solution methods.

We contribute to \colr{the research field of appointment scheduling} in \colr{five} ways. 
First, we build on a model for timeslot allocation, presented by \citet{Hulshof2013, Hulshof2016} \colr{and extend the contribution function with rewards for appointments and surgeries to stimulate the optimal use of the available capacity. 
Second, we compare the performance of multiple solution methods. We evaluate four static decision rules, }formulate a similar (I)LP to \citet{Hulshof2013}, and investigate a form of approximate dynamic programming suitable for infinite-time-horizon problems: least-squares policy iteration (LSPI). These methods are implemented on data from SMK and their performance is compared with the current planning method used by the hospital.
\colr{Third}, we implement the concept of planning ahead. \colr{To accommodate timely booking of appointments, timeslot} schedules are generated weeks in advance. Perfect information about the future state of the system for which we wish to make a planning decision is then unavailable. In this paper, we develop a method to make a prediction of that future state. Furthermore, we compare the performance of the introduced solution methods when this planning horizon is increased. 
\colr{Fourth}, to the best of our knowledge, the variant of LSPI used in this paper has not been used previously in literature. An off-policy variant using an approximate Q-function is introduced in \citet{Lagoudakis2003}, and an on-policy variant using an approximate value function is presented in \citet{Powell2019}. We introduce and implement an off-policy variant using an approximate value function. 
\colr{Fifth}, rather than using a post-decision state within the algorithm as is done in \citet{Powell2019} and \citet{Hulshof2016}, we argue why an expectation of the next state can be used. We also extensively investigate the influence of adaptations to hyperparameters and to the algorithm itself on convergence and performance on the timeslot allocation problem.

\colr{In Section \ref{sec:related_research}, we provide an overview of related research in the field of healthcare planning. Section \ref{sec:model} presents our model for the timeslot allocation problem. Section \ref{sec:solution_methods} introduces the solution methods which are used. Section \ref{sec:data_analysis} describes the small test instance, large test instance, and the SMK instance we use to test the solution methods. We discuss the results of the solution methods per instance in Section \ref{sec:results} and present our conclusions and recommendations in Section \ref{sec:conclusion_recommendations}. For a list of abbreviations, see Appendix \ref{sec:abb}.}

%% file: literature.tex
\section{Literature review}\label{sec:related_research}
In this section, we discuss relevant research related to outpatient appointment scheduling and the timeslot allocation problem. Section \ref{sec:outpatient} discusses literature on outpatient appointment scheduling. Section \ref{sec:od_or_thesis} presents previous research at the SMK on planning OD and OR sessions. 
Section \ref{sec:hulshof} describes research on the timeslot allocation problem, and two solution methods, on which our \colr{approach} is based.

\subsection{Outpatient appointment scheduling}\label{sec:outpatient}

Outpatient appointment scheduling has been intensively studied, see for an overview \citet{Ahmadi2017} and \citet{Ala2022}. In this section, we discuss the research related to our timeslot allocation problem. 

\citet{Erdelyi2011} consider the case where there is a fixed amount of daily processing capacity. Each day, jobs with different priorities arrive randomly for which it needs to be decided which jobs are scheduled on which days. The waiting cost of a job depends on its priority level. The objective is to minimize the total expected waiting cost over a finite planning horizon. This setting differs from our setting, because jobs arriving on different days are considered to be independent. In our case, patients need follow-up appointments, which means that the jobs are not necessarily independent. Therefore, in \citet{Erdelyi2011}, no transition probabilities between states needs to be taken into account. The problem is solved using approximate dynamic programming where the value function is approximated by decomposing the problem per day. \citet{Patrick2008} consider a similar problem for which approximate dynamic programming techniques are used to formulate an approximate linear program, which is solved using column generation. The resulting policy is tested through simulation.

\citet{Nunes2009} consider a similar problem to \colr{timeslot allocation} in which patients from different specialties are admitted during fixed planning periods. Each patient has a treatment pattern using resources during their admission. The problem is modelled as a Markov decision process for which the action is the number of patients to admit for each speciality in the next planning period. The state is the number of patients for each specialty that has a specific treatment pattern during the last period. The transition probabilities indicate the probability that a patient gets a different treatment pattern. The goal is to be as close as possible to the target utilization of the resources. This is modelled through a cost function consisting of costs for underutilization (under target), overutilization (above target), and overcapacity (demand exceeds capacity). Only small instances are solved using value iteration, and the development of solution methods for larger instances is left as future work.

\citet{Nguyen2015} determine the capacity required to schedule first and follow-up appointments while respecting access time targets. For both types of appointments, a constant discharge rate is considered. This problem differs from the setting considered in our paper, as we consider the capacity to be given and the access time to be minimized. \citet{Aslani2021} extend upon the work of \citet{Nguyen2015} by considering the demand for first appointments to be uncertain. They develop a robust optimization model that minimizes the maximum capacity needed for the worst case realization of demand.

\citet{Bakker2017} dynamically allocate capacity to multiple specialists and activities. The activities consist of first appointments, follow-up appointments, and surgery. The capacity is \colr{not assigned optimally, but approximated} using a greedy heuristic aimed at minimizing the variability in roster activities. The approach is evaluated using simulation.

\citet{Laan2018} optimise appointment scheduling with respect to access time, taking fluctuating patient arrivals and unavailabilities of physicians into account. Follow-up appointments are modelled as a patient type, which means that these appointments are modelled independently, while this is not the case in practice. They formulate a stochastic mixed integer programming problem, and approximate its solution using two different approaches: (1) a mixed integer programming approach that results in a static appointment schedule, and (2) Markov decision theory, which results in a dynamic scheduling strategy. They apply the methodologies to a case study of the surgical outpatient clinic of the Jeroen Bosch Hospital. 

\citet{Deglise2018} develop an appointment template for an integrated care environment where multiple specialties serve multiple types of patients. In this template, new patients are scheduled and capacity for follow-up downstream appointments is \colr{not planned but} reserved based on a stochastic location function. The problem is formulated as a queuing network optimization \colr{problem} and an approximate solution is determined by solving it as a deterministic linear optimization problem.

\subsection{OD and OR session allocation} \label{sec:od_or_thesis}
\cite{Hattingh2019} and \cite{Tsai2017} \colr{discuss work} at the Sint Maartenskliniek to determine the optimal allocation of OD and OR sessions for orthopaedic surgeons. They aim to stabilize access times for patients requiring surgery, reduce unused OR time, and meet the yearly production targets of each surgeon. Access times were stabilized by reducing fluctuations around a target OR waiting list length. Both modelled the patient pathway at SMK as a Markov decision process, but applied different methods to solve the problem.

\citet{Tsai2017} developed a detailed model which included the diagnostics and screening process. Stochastic dynamic programming was applied to determine an optimal allocation, but this was impractical due to high computation times for calculating transition probabilities between states. Even after \colr{signicantly reducing the state space through}  simplifying assumptions, the model remained computationally expensive. \citet{Tsai2017} recommends ignoring the radiology department and mentions approximate dynamic programming as a possible solution method.

\citet{Hattingh2019} formulated a similar MDP model to Tsai, following the recommendation to ignore the radiology department. Test cases where the effects of OD sessions and OR sessions on waiting lists were modelled as either stochastic or deterministic were compared. Replacing the stochastic parameters by deterministic averages resulted in similar solutions for the test cases and significantly reduced the running time. This assumption allowed for the formulation of an integer linear program in order to optimize the model. The ILP was successfully able to determine the number of OD and OR sessions per week for each surgeon in a two-week planning period. \colr{These outcomes are input for our timeslot allocation problem in which the given number of OD and OR sessions are assigned to the different appointment types.}

\subsection{Patient admission planning} \label{sec:hulshof}
A very similar problem to the one under consideration in this paper was investigated by \citet{Hulshof2013} and \citet{Hulshof2016}. They describe a general model for tactical resource allocation and patient admission planning in health care, which is formulated as a \colr{mixed integer linear program (MILP) in \citet{Hulshof2013} and as an Markov Decision Process (MDP)} in \citet{Hulshof2016}. 

\citet{Hulshof2013} formulated the problem as a mixed integer linear program. Since an MILP cannot solve stochastic problems, the transition probabilities are modelled as deterministic fractions $q_{ij}$ of patients moving from one stage (or queue) to another. The objective is ``to achieve equitable access and treatment duration for patient groups and to serve the strategically agreed number of patients''. In the objective function, this is approximated by minimising the weighted sum of patients waiting in each queue over all access times and time periods.
The authors apply an iterative scheme in order to determine the weights used in the objective function, so that the true performance targets set by the hospital are optimized. For relatively large instances (fifty queues, two resources, and six time periods), calculation time took on average four minutes with an integrality gap of $0.01\%$. Results show that the model divides resources more equitably over time and the number of patients served is closer to the strategically set target. The authors recommend implementing this model using a rolling-horizon approach, i.e., only applying decisions for earlier time periods while considering expected effects on later time periods. The model is relatively easy to adapt with, for example, more constraints or a different objective function.

\citet{Hulshof2016} formulated the problem as an MDP. The objective is to generate a schedule which provides equitable access to treatment for patients and meets the strategic goals of the hospital. 
Computing the exact solution to this MDP using dynamic programming is only possible for very small instances, as the state, decision and outcome spaces grow exponentially in size (the three curses of dimensionality). Therefore, approximate dynamic programming (ADP) is applied, which is a form of approximate model-free on-policy value iteration. The algorithm overcomes intractability by using a post-decision state and a linear parametrisation of the value function using basis functions to approximate future costs. An ILP is used to efficiently find an optimal action given the current approximate value function and problem constraints. Various basis functions are proposed and compared using regression analysis. Recursive least squares is used to update the parameter vector used for the value function.
The authors use small test instances to show that the ADP algorithm provides a close approximation to the outcome of exact dynamic programming in reasonable time. For large instances (forty queues, four resources, eight time periods), the algorithm takes approximately one hour to converge to a parameter vector. The resulting policy performs much better on test instances than two proposed greedy heuristics. The authors conclude that ADP is an appropriate technique for tactical healthcare planning.

%% file: model.tex
\section{Problem formulation}\label{sec:model}
\colr{In this section, we describe the model for the timeslot allocation problem. Section \ref{sec:problem_definition} describes the problem and discusses our assumptions. Section \ref{sec:base_model} describes the sets, parameters, functions, states and actions defining the model. Section \ref{sec:model_objective} proposes the objective function and Section \ref{sec:model_planning_ahead} describes how we incorporate the concept of planning ahead.
All notation used in this section is summarized in Table \ref{tab:model_svp}.}

\subsection{Problem definition}\label{sec:problem_definition}

\colr{We aim to model the effects of allocating timeslots to different appointment types on the waiting lists of a single surgeon. The input is the number of OD and OR timeslots available in the planning time period, which can be determined from the activity plan generated in the tactical planning phase. Relevant parameters for the considered surgeon, and the current state of the waiting lists, should be given as well. The model should output how many OD and OR timeslots should be allocated to which appointment type per time period. 

Since treating patients not only has an effect on the current waiting lists but also on the future, as patients move through their treatment process, future consequences of current decisions must be taken into account. Additionally, since no patient is the same, these consequences are uncertain and we are faced with a sequential decision problem under uncertainty. For this reason, the problem is modelled as a Markov decision process (MDP). 

We use the }following assumptions:

\begin{enumerate}
	\item The total number of available OD and OR timeslots are given for each time period. There is sufficient capacity of other resources (locations, materials) to utilize all timeslots.
	\item The duration of an OD and OR timeslot is fixed per surgeon. We allow for OD appointments to take up multiple OD timeslots. For OD appointments, this is a reasonable assumption. This is not necessarily the case for OR appointments, as the time required depends on the complexity of the surgery. We use the rounded average of surgery durations per specific surgeon.
	\item We ignore diagnostics in the radiology department and appointments for pre-operative screening stages.
	\item The considered surgeon is not affected by, or dependent on, activities of other surgeons. \colr{In the SMK, the majority of patients stay with the same surgeon for all consultations and surgery.}
	\item Patient appointments are booked into appointment slots in the most efficient way. \colr{This means that we assume that the patient resulting in the highest contribution is always available and booked. In real life, the highest contribution patient may not be available at the required time, appointments can be cancelled, and patients with higher urgency may arrive later on, when there are no available timeslots. However, to compare the quality of solution methods, this assumption suffices.}
        \item We assume that all appointments are realised (zero no-shows), \colr{which also suffices when we use the model to compare the quality of solution methods.}
\end{enumerate}

\subsection{Model}\label{sec:base_model}

Our model is based on the model by \citet{Hulshof2016}. \colr{Different from \citet{Hulshof2016}}, we consider an infinite time horizon $\mathcal{T}$ in order to take into account all future consequences of actions. Unless stated otherwise, we set $\mathcal{T} = \mathbb{Z}^+$. The patient queues are \colr{given by set $\mathcal{J}$. Another addition to the work of \citet{Hulshof2016}, is that each queue $j\in \mathcal{J}$ has a given access time target denoted by $u_j$}. The number of time periods that a patient has been waiting in a queue is indicated by the value $w \in \mathbb{Z}^+$. For some solution methods, it may be necessary to upper bound the \colr{waiting} time \colr{for queue $j\in\mathcal{J}$} with some value $W_j$. \colr{Therefore, we introduce the set $\mathcal{W}_j$, which is equal to either $\mathbb{Z}^+$ or $\lbrace 0,...,W_j\rbrace$ depending on the solution method.}
The number of new patients arriving from outside the system to queue $j \in \mathcal{J}$ at time $t \in \mathcal{T}$ is given by $\lambda_{j, t}$. A patient in queue $i \in \mathcal{J}$ transfers to queue $j \in \mathcal{J}$ with probability $q_{i, j}$ when treated. \colr{A patient in queue $i\in \mathcal{J}$ leaves the system with probability $1-\sum_{j\in \mathcal{J}}q_{i,j}$.}  \colr{When a patient is not treated, the waiting time $w\in \mathcal{W}_j$ is increased by one.}  We consider a set of resources $\mathcal{R} = \{$OD, OR$\}$ \colr{consisting of} the outpatient department and operation room. The resource capacity $\eta_{r, t}$ indicates the number of timeslots available for resource $r \in \mathcal{R}$ at time $t \in \mathcal{T}$. \colr{Parameter} $\zeta_{j, r}$ indicates how many timeslots of resource $r \in \mathcal{R}$ are required by a patient in queue $j \in \mathcal{J}$. 

The sub-state $s_{j, w, t}$ represents the number of patients waiting in queue $j \in \mathcal{J}$ for \colr{$w\in \mathcal{W}_j$} time periods at time $t \in \mathcal{T}$. The state $s_t$ is a vector $(s_{j, w, t})_{j \in \mathcal{J}, w \in \mathcal{W}_j}$. 
The size of a state $\abs{s_t}$ is the sum of its sub-states, and represents the total number of patients in the system at time $t$. The state space of all possible states is denoted by $\mathcal{S}$.
The sub-action $a_{j, w, t} \in \mathbb{Z}^+$ determines how many patients from sub-state $s_{j, w, t}$ to treat. Action $a_t$ is a vector $(a_{j, w, t})_{j \in \mathcal{J}, w\in\mathcal{W}_j}$. This is the decision variable of the problem. \colr{Note that we have as restriction that we cannot treat more patients than are available in queue $j\in \mathcal{J}$, waiting for $w\in \mathcal{W}_j$ time periods at time  $t\in \mathcal{T}$. For each resource $r\in \mathcal{R}$, we also make sure that the needed capacity to treat patients does not exceed the  available capacity  at time $t\in \mathcal{T}$.} The action space of all \colr{feasible} actions is denoted by $\mathcal{A}$. The set of \colr{feasible} actions given a state $s_t$ \colr{at time $t\in \mathcal{T}$} is denoted by $\mathcal{A}(s_t)$ \colr{and is given by:

\begin{equation}
\mathcal{A}(s_t)=\left\lbrace a_t : \begin{array}{ll}
a_{j,w,t}\leq s_{j,w,t},    & \forall j\in\mathcal{J},w\in \mathcal{W}_j\\
\displaystyle\sum_{j\in\mathcal{J}}\sum_{w\in\mathcal{W}_j}\zeta_{j,r}a_{j,w,t}\leq \eta_{r,t},    &\forall r\in\mathcal{R} 
\end{array}
\right\rbrace
\label{eq:feasibleactions}
\end{equation}}

\colr{Given a state $s_t$ and action $a_t$ at time $t\in \mathcal{T}$, the resulting next state $s_{t+1}$ is given by:

\begin{equation}
s_{j,w,t+1}=\left\lbrace
\begin{array}{ll}
\lambda_{j,t+1}+\displaystyle\sum_{i\in\mathcal{J}}\sum_{v\in \mathcal{W}_j}q_{i,j}a_{i,v,t},&w=0,\forall j\in \mathcal{J}\\
s_{j,w-1,t}-a_{j,w-1,t},&1\leq w\leq W_j-1,\forall j\in \mathcal{J} \\
\displaystyle\sum_{v=W_j-1}^{W_j}s_{j,w,t}-a_{j,w,t},&w=W_j, \forall j\in \mathcal{J}
\end{array}\right.
\label{eq:nextstate}
\end{equation}}

\begin{table}
\caption{\label{tab:model_svp}\colr{Sets, variables, parameters and functions used for the model}}
\begin{tabular}{ll}
\hline
\textbf{Sets}       &  \\
$\mathcal{T}$       & Time periods  \\
$\mathcal{J}$       & Queues\\
\colr{$\mathcal{W}_j$}       & \colr{Waiting times for queue $j\in \mathcal{J}$}\\
$\mathcal{R}$       & Resource types  \\
$\mathcal{S}$       & State space  \\
$\mathcal{A}$       & Action space   \\
                    &    \\
\textbf{Variables}  &     \\
$a_{j, w, t}$      & \begin{tabular}[c]{@{}l@{}}Sub-action: number of patients to treat at time $t \in \mathcal{T}$ from queue $j \in \mathcal{J}$ who have been waiting \\for \colr{$w\in\mathcal{W}_j$} time periods\end{tabular} \\
$s_{j, w, t}$      & \begin{tabular}[c]{@{}l@{}}Sub-state: number of patients in queue $j \in \mathcal{J}$ who have been waiting for \colr{$w\in\mathcal{W}_j$} time periods \\at time $t \in \mathcal{T}$\end{tabular}                        \\  &  \\
\textbf{Parameters} &  \\
\colr{$u_j$} & \colr{Access time target of patients in queue $j\in \mathcal{J}$}\\
$\eta_{r,t}$        & Available \colr{number of timeslots for} resource $r \in \mathcal{R}$ at time $t \in \mathcal{T}$  \\
$\zeta_{j,r}$           & \colr{Number of timeslots} of resource $r \in \mathcal{R}$ required by patient in queue $j \in \mathcal{J}$    \\
$q_{i, j}$      & \begin{tabular}[c]{@{}l@{}}Probability that a patient transfers from queue $i \in \mathcal{J}$ to queue $j \in \mathcal{J}$\end{tabular}  \\
$\lambda_{j, t}$  & Number of new patients to enter queue $j \in \mathcal{J}$ at time $t \in \mathcal{T}$   \\
\colr{$\omega_j$}&\colr{Relative weight of queue $j\in \mathcal{J}$}\\
\colr{$c_{j,w}$}&\colr{Cost of patients in queue $j\in \mathcal{J}$ who have been waiting for $w\in\mathcal{W}_j$ time periods}\\
\colr{$r_j$}&\colr{Reward for treating a patient from queue $j\in \mathcal{J}$}\\
\colr{$\gamma$} & \colr{Discount factor}\\
\colr{$p$}        & Planning horizon: the number of time periods to plan ahead        
                      \\  &  \\
\colr{\textbf{Functions}} &  \\
\colr{$C(s_t,a_t)$} & \colr{Contribution function for taking action $a_t$ in state $s_t$}\\
\colr{$\pi(s_t)$} & \colr{Policy function defining which action to take in state $s_t$}\\
\colr{$V(s)$} & \colr{Value function of state $s$}\\
 \hline
\end{tabular}
\end{table}

\subsection{Objective} \label{sec:model_objective}
\colr{The objective of our MDP model is twofold: maximize the number of patients treated while minimizing the access time of patients treated after their access time target}. The objective function \colr{is a}combination of these two \colr{aspects}.

The \colr{first} objective\colr{, in which our work differs from the work of \citet{Hulshof2016},} is to maximize the number of patients treated \colr{as this increases both the patients' health benefits and resource utilization}.  Therefore, a reward $r_{j}$ is assigned for every patient treated from queue $j \in \mathcal{J}$. As a result, the \colr{first} objective is given by:
\begin{equation*}
\max \sum_{t \in \mathcal{T}} \sum_{j \in \mathcal{J}} \sum_{w \in \mathcal{W}_j} r_{j} \cdot a_{j, w, t}.
\end{equation*}

\colr{The second objective can be expressed as minimizing the weighted sum of patients who are not treated on time. We define \colr{the cost of letting a patient in queue $j\in \mathcal{J}$ wait for $w\in\mathcal{W}_j$ time periods by $c_{j,w}$}. If the waiting time $w\in\mathcal{W}_j$ has not yet exceeded the access time target $u_j$, the patient from queue $j\in\mathcal{J}$ is not late yet and the cost should be zero: $c_{j, w} = 0$ for $w< u_j$.
Furthermore, $c_{j, w}$ should be non-decreasing in $w\in\mathcal{W}_j$ (as higher access times should be penalized more heavily), and non-increasing in $u_j$ (as a higher value of $u_j$ indicates lower urgency). A possible cost function for a patient in queue $j\in\mathcal{J}$ waiting for $w\in\mathcal{W}_j$ could be $c_{j, w} = \omega_j \cdot \frac{w}{u_j}$ for $w \geq u_j$, to reflect the waiting time relative to the access time target. Here, $\omega_j$ is a parameter determining the weight of each queue relative to the other queues.
Now, the second objective is given by:
\begin{equation*}
\min \sum_{t \in \mathcal{T}} \sum_{j \in \mathcal{J}} \sum_{w\in \mathcal{W}_j} c_{j, w} \cdot (s_{j, w, t} - a_{j, w, t}).
\end{equation*}}

\colr{This means that each patient who has been waiting in queue $j\in \mathcal{J}$ for $w\in\mathcal{W}_j$ time periods at time $t\in \mathcal{T}$ and is not chosen to be treated incurs a cost of $c_{j,w}$. Since these patients transition to state $s_{j,\min(W_j,w+1),t+1}$, they might also incur a cost of $c_{j,\min(W_j,w+1)}$ in the next time period. So a patient from queue $j\in \mathcal{J}$ who is chosen to be treated after waiting for $w\in\mathcal{W}_j$ time periods incurs a total cost of $\sum_{v\in \mathcal{W}_j: v \geq u_j}^{w} c_{j,v}$.}

Combining the two objectives gives the following contribution function \colr{that we aim to maximize}:
\begin{equation} \label{eq:contribution_function}
C(s_t, a_t) = \sum_{j \in \mathcal{J}}  \sum_{w \in \mathcal{W}_j} r_{j} \cdot a_{j, w, t} - c_{j, w} \cdot (s_{j, w, t} - a_{j, w, t}).
\end{equation}
Note that if the reward is set too low in comparison to the cost of waiting patients, then the optimal solution could be to accept very few new patients in order to enforce that patients never wait too long. On the other hand, if it is set too high, it will become beneficial to have  very long waiting lists so that there are always enough patients to fill available resource capacity. This trade-off must be considered carefully.

The \colr{aim} is, given initial state $s \in \mathcal{S}$, to find an optimal policy $\pi^*$ maximizing the value function. \colr{Here, policy $\pi$ defines for each state $s_t$ which action $a_t$ to take, i.e., $\pi(s_t)=a_t$. The discount factor is given by $\gamma$.}
\begin{equation*}
V^*(s) = \max_{\pi}\mathbb{E}_{\pi} \left[ \sum_{t \in \mathcal{T}} \gamma^t C(s_t, \pi(s_t)) | s_0 = s \right] 
\end{equation*}

\subsection{Planning ahead} \label{sec:model_planning_ahead}
Until this point, we have assumed that an action \colr{$a_t$} is selected at time \colr{$t\in \mathcal{T}$} with perfect knowledge of the state \colr{$s_t$} at time \colr{$t\in\mathcal{T}$}. However, planning decisions are made in advance, so that surgeons are informed of their schedules, and patients of their appointments, on time. The number of time periods to plan ahead is defined as the \textit{planning horizon} and is denoted by the parameter $p$. This implies that at time \colr{$t\in \mathcal{T}$}, we must determine action $a_{t + p}$ \colr{based on incomplete knowledge of state $s_{t+p}$}. To \colr{obtain an estimate} of state $s_{t+p}$\colr{, we can use} the current state $s_t$, and the actions $a_t, \ldots, a_{t+p-1}$ that \colr{will be} taken in the meantime, as these actions have \colr{already} been chosen in previous time periods. \colr{When introducing the various solution methods in the next section, we do not specifically take planning horizon $p$ into account. Planning ahead can be applied to any of the introduced solution methods by not deciding on an action for the next time period, but for $p$ time periods ahead. The effect of planning ahead on the solution quality is discussed in Section \ref{sec:results_smk}.}

%% file: solution_methods.tex
\newpage
\section{Solution methods}\label{sec:solution_methods}
\colr{As} we have modelled the problem as a Markov decision process, the solution takes the form of a policy: a function or rule dictating which action to take \colr{for each} state. An optimal policy maximizes the expected contribution gained by following that policy over time. Such an optimal policy can be found through exact value iteration (EVI), which we discuss in Section \ref{sec:solution_evi}. However, \colr{as the state, decision, and outcome spaces grow exponentially in size (the three curses of dimensionality),} we investigate \colr{several} approximation methods and heuristics. These are least-squares policy iteration (Section \ref{sec:solution_lspi}), integer linear programming (Section \ref{sec:solution_lp}), and decision rules (Section \ref{sec:solution_decisionrule}). \colr{Finally, Section \ref{sec:solution_static} describes the current (static) allocation method used by SMK that we use as a benchmark.}

\subsection{Exact value iteration}\label{sec:solution_evi}
The exact solution to an MDP can be found using dynamic programming. For realistic instances of the timeslot allocation problem, this will not be possible, as the three curses of dimensionality all apply. However, it is useful to find the exact solution to small test instances, to compare the performance of the other \colr{proposed} solution methods. We also use this to compare different basis functions that can be used for least-squares policy iteration.

Finding an exact solution requires solving the Bellman equation:
\begin{equation}
V^*(s) = \max_{a \in \mathcal{A}(s)} \left( C(s, a) + \gamma \sum_{s' \in \mathcal{S}} P(s'|s, a) V^*(s') \right) \text{, } \forall s \in \mathcal{S}.\label{eq:val_iteration_bellman}
\end{equation}

This \colr{equation can be solved} using the exact value iteration algorithm for infinite horizon optimization \colr{introduced by} \cite{Powell2019}. However, for our problem, it is not \colr{straightforward} to calculate \colr{probability $P(s'|s, a)$ to transition from state $s$ to state $s'$ under action $a$}. An action causes patients to leave and enter queues with given probabilities, but calculating the probability that all these transitions add up to a \colr{certain} state \colr{$s'$} turns out to be complex.
A \colr{computationally less intensive approach} approach, giving an equivalent solution, is not to sum over all possible next states, but over all possible realizations of transitions caused by action $a$. Let $I = \{(i, j): q_{ij} > 0, i, j \in \mathcal{J} \cup \{0\} \}$ be the pairs of indices of queues between which patients can transition, including the `\colr{exit}' queue $0$. Now let $X = (X_{i, j})_{(i, j) \in I}$ be a random variable denoting the number of patients transitioning from queue $i$ to queue $j$. Let $a_i$ be the action indicating the number of patients to treat from queue $i \in \mathcal{J}$. Then, $X_{ij}$ is binomially distributed: $X_{ij} \sim B(a_i, q_{ij})$. \colr{Now}, let $\mathcal{X}(a)$ be the set of possible realizations of transitions resulting from action $a$\colr{, i.e.,} 
\begin{equation*}
\mathcal{X}(a) = \left\{ x \in \mathbb{Z}^{\abs{I}} : x_{ij} \geq 0, \text{ } \forall (i, j) \in I \text{ and } \sum_{j \in \mathcal{J}} x_{ij} = a_i, \text{ } \forall i \in \mathcal{J} \right\}.
\end{equation*}
Using the transition function $T(s, a, x)$\colr{, which} provides the new state given previous state $s$, action $a$, and (a realization of) exogenous information $x$, Equation (\ref{eq:val_iteration_bellman}) can now be written as:
\begin{align*}
V^*(s) &= \max_{a \in \mathcal{A}(s)} \left( C(s, a) + \gamma \sum_{x \in \mathcal{X}(a)} P(X = x | a) V^*(T(s, a, x)) \right) \\
&= \max_{a \in \mathcal{A}(s)} \left( C(s, a) + \gamma \sum_{x \in \mathcal{X}(a)} \left( \prod_{(i,j) \in I} P(X_{ij} = x_{ij} | a_i) \right) V^*\left(T(s, a, x)\right) \right) \\
&= \max_{a \in \mathcal{A}(s)} \left( C(s, a) + \gamma \sum_{x \in \mathcal{X}(a)} \left( \prod_{(i,j) \in I} \binom{a_i}{x_{ij}} q_{ij}^{x_{ij}} (1 - q_{ij})^{a_i - x_{ij}} \right) V^*\left(T(s, a, x)\right) \right),
\end{align*}

\colr{\noindent which can be solved using the exact value iteration algorithm.}

\subsection{Least-squares policy iteration}\label{sec:solution_lspi}

We consider least-squares policy iteration as an approximation method \colr{to solve larger instances}. The algorithm \colr{determines} a policy by learning a parametrised approximation of the value function, $\overline{V}(s) = \theta^{\top}\phi(s)$, with $\phi(s)$ a \colr{basis function} vector and $\theta$ a parameter vector. LSPI is hindered less by the curses of dimensionality, and can therefore be used to solve larger instances of the timeslot allocation problem.

\colr{We first give the general framework of our LSPI implementation in Section \ref{sec:general_framework}. In Section \ref{sec:simulation_environment}, we describe how we sample the states and exogenous information. To determine the next action, we solve the linear program given in Section \ref{sec:action_selection}. The basis functions are described in Section  \ref{sec:basis_functions}. }

\subsubsection{\colr{General framework}}\label{sec:general_framework}
We propose a combination of the off-policy LSPI-Q \citep{Lagoudakis2003} and on-policy LSPI-V \citep{Powell2019} algorithms, \colr{resulting in} the off-policy LSPI-V algorithm shown in Algorithm \ref{alg:off_lspi_v} where the approximate transition function $\overline{T}(s, a)$ is used as an approximation of the next state. To the best of our knowledge, this version has not been used yet in the literature. The value function version must be used as the action space is \colr{too} large \colr{to approximate} the $Q$-function. The off-policy variant is used, because the on-policy variant did not converge. Finally, we use the least-squares temporal differencing (LSTD) update \citep{Bradtke1996}, as opposed to  recursive LSTD, \colr{to reduce the number of} matrix-vector multiplications. The number of basis functions $F$ for this problem is small enough to efficiently compute the inverse of $A_M$, with $M$ the number of state samples. 

\begin{algorithm}[h]
\SetAlgoLined
\LinesNumbered
$\theta_0 = \mathbf{0} \in \mathbb{R}^F$\colr{, $n=0$} \\
Choose small constants $\epsilon, \delta > 0$ \\
Generate a \colr{set} of \colr{sampled} states $D = \{s_1, \ldots, s_M\}$\\
Define the policy:
$\pi\left(s | \theta \right) = \arg\max_{a \in \mathcal{A}(s)} \left( C(s, a) + \gamma \theta^{\top} \phi(\overline{T}(s, a)) \right)$\\
\While{$|| \theta_n - \theta_{n-1}|| \geq \delta$\colr{ or $n=0$}}{
	\colr{$n=n+1$}\\
 $b_0 = \mathbf{0} \in \mathbb{R}^F$ \\
	$A_0 = \epsilon I$, using $F \times F$ identity matrix $I$ \\
	\For{$m = 1, \ldots, M$}{
		Compute action $a_{m} = \pi(s_{m} | \theta_{n-1})$ \\
		Obtain a sample of exogenous information $x_m$ from the simulation \\
		Determine the next state $s'_{m} = T(s_{m}, a_{m}, x_m)$ \\
		Perform the LSTD update:
		$b_m$ = $b_{m-1} + C(s_m, a_m) \phi(s_m),
		A_m = A_{m-1} + \phi(s_m) \cdot (\phi(s_m) - \gamma \phi(s'_m))^{\top}$\\
	}
	$\theta_n = A_{M}^{-1} b_{M}$ \\
}
$N=n$\\
\KwRet{the approximate policy $\pi\left(s|\theta_N \right)$}
\caption{Off-policy LSPI-V}\label{alg:off_lspi_v} 
\end{algorithm}

\subsubsection{\colr{Sampling of states and exogenous information}}\label{sec:simulation_environment}
\colr{LSPI requires generating a set of $M$ sampled states $D=\lbrace s_1,\ldots,s_M\rbrace\subseteq \mathcal{S}$. A sampled state is obtained by first randomly generating the total number of patients in the system. For each of these patients in the system, a random patient care pathway is generated. A patient care pathway contains the ordered appointments/queues of a patient. Next to this, we generate the patient's current queue in their pathway and current waiting time.  Given the current queue and waiting time of these patients, we can determine the number of patients in each queue $j\in \mathcal{J}$ waiting for $w\in \mathcal{W}_j$ time periods giving us a sampled state.

In addition, to determine the next state $s'_m$, we need to sample exogenous information $x_m$, which contains the number of newly arriving patients for each queue and to which queue treated patients transition. When a patient is treated, the next queue can easily be determined from the generated care pathways for this sampled state $s_m$.}

\subsubsection{Action selection}\label{sec:action_selection}
Using the policy function to choose an action requires solving:
\begin{equation*}
\pi\left(s | \theta \right) = \arg\max_{a \in \mathcal{A}(s)} \left( C(s, a) + \gamma \theta^{\top} \phi(\overline{T}(s, a)) \right),
\end{equation*}
where $\overline{T}(s, a)$ is an approximation of the next state. 
When an MDP can be solved exactly, given optimal value function $V^*$, the optimal policy can be described as:
\begin{align*}
\pi(s) &= \arg \max_{a \in \mathcal{A}(s)} \left( C(s, a) + \gamma \sum_{s' \in \mathcal{S}} P(s'|s, a)V^*(s') \right) \\
&= \arg \max_{a \in \mathcal{A}(s)} \left( C(s, a) + \gamma \mathbb{E}[V^*(s') | s, a] \right).
\end{align*}
Now, rather than an optimal value function, we have an approximate value function $\overline{V}$. Since we are using a linear parametrisation and linear basis functions, this function is linear in the state vector. Therefore, we can rewrite $\mathbb{E}[\overline{V}(s') | s, a]$ \colr{to} $\overline{V}(\mathbb{E}[s' | s, a])$. This gives:
\begin{equation*}
{\pi}\left(s|\theta\right) = \arg \max_{a \in \mathcal{A}(s)} \left( C(s, a) + \gamma \theta^{\top}\phi(\mathbb{E}[s' | s, a]) \right).
\end{equation*}

For this problem, the expected next state can be calculated \colr{using the transition probabilities $q_{i,j}$ from queue $i\in\mathcal{J}$ to queue $j\in \mathcal{J}$ and expected number of newly arriving patients $\lambda_j$ to queue $j\in\mathcal{J}$. }  Given an action $a$, the expected number of patients transitioning from queue $i\in\mathcal{J}$ to queue $j\in\mathcal{J}$ is \colr{given by} $\sum_{w\in\mathcal{W}_j} q_{i, j} \cdot a_{i, w}$. The rest of the state transition\colr{s} is determined by deterministic events \colr{as patients who are not treated} wait \colr{an} extra time period.

\colr{Using the above, we can formulate the following integer linear program to determine the next action resulting in  the expected next state $\overline{s}$, which we refer to as the \textit{policy-LP}:}
\begin{maxi!}<b>
{a}{ C(s, a) + \gamma \sum_{f = 1}^F \theta_f \cdot \phi_f(\overline{s}) }{\label{lp:policy}}{ }
\addConstraint{\overline{s}_{j, 0}}{ = \lambda_{j} + \sum_{i \in \mathcal{J}} \sum_{w\in\mathcal{W}_i} q_{i, j} \cdot a_{i, w} , \label{constr:policy_transition_prob}}{\quad \forall j \in \mathcal{J}}
\addConstraint{\overline{s}_{j, w}}{ = s_{j, w-1} - a_{j, w-1} , \label{constr:policy_treat}}{\quad \forall j \in \mathcal{J}, 1 \leq w \leq W_j - 1}
\addConstraint{\overline{s}_{j, W_j}}{ = \sum_{w = W_j-1}^{W_j} s_{j, w} - a_{j, w} , \label{constr:policy_treat_W}}{\quad \forall j \in \mathcal{J}}
\addConstraint{a_{j, w}}{ \leq s_{j, w}, \label{constr:policy_action_queuelength}}{\quad \forall j \in \mathcal{J}, w\in \mathcal{W}_j}
\addConstraint{\sum_{j \in \mathcal{J}}  \sum_{w\in\mathcal{W}_j} \zeta_{j, r}a_{j, w}}{\leq \eta_r , \label{constr:policy_action_resources}}{\quad \forall r \in \mathcal{R}}
\addConstraint{a_{j, w}}{ \geq \colr{0}, \label{constr:policy_action_integer}}{\quad \forall j \in \mathcal{J}, w\in\mathcal{W}_j}
\addConstraint{\colr{\overline{s}_{j, w}}}{ \geq \colr{0},}{\colr{\quad \forall j \in \mathcal{J}, w\in\mathcal{W}_j.}}
\end{maxi!}
The contribution function $C(s,a)$ is as defined in Equation \eqref{eq:contribution_function}. Constraint\colr{s} \eqref{constr:policy_transition_prob}- \eqref{constr:policy_treat_W} \colr{determine the expected next state following Equation \eqref{eq:nextstate}}. Constraint\colr{s} \eqref{constr:policy_action_queuelength}-\eqref{constr:policy_action_resources} \colr{define all feasible actions following Equation \eqref{eq:feasibleactions}. Note that for $\lambda_j$ and $\eta_r$, we take the expected value of $\lambda_{j,t}$ and $\eta_{r,t}$ over $t\in\mathcal{T}$, respectively.}

\colr{Even though it is only possible to treat an integer number of patients, we have relaxed the restriction on action variable $a_{j,w}$.} Since this program must be solved for every state in dataset \colr{$D$}, in every iteration of the algorithm, we \colr{have relaxed} this integer constraint to a non-negativity constraint in order to reduce the solving time. \colr{This means that} the actions resulting in the optimal solution should be rounded down to the nearest integer to ensure feasibility \colr{as} rounding up may violate constraints \eqref{constr:policy_action_queuelength} or \eqref{constr:policy_action_resources}. Rounding down could result in loss of optimality, as some available capacity could go unused. When the policy is actually used for planning, the integer constraint can and should be used, as in that case the policy-LP only needs to be solved once per planning period.

\subsubsection{Basis functions}\label{sec:basis_functions}
Well-chosen basis functions are essential \colr{for} the convergence and performance of LSPI. In this section, various basis functions are proposed. We also introduce a linear regression method to \colr{determine which basis function best approximates the exact value function $V^*$.}

Basis functions may be non-linear (for example, Gaussian or quadratic functions of features), but we focus on linear functions to ensure that an action selection can be done by solving the policy-LP, and such that $\mathbb{E}\left[\phi(s)\right] = \phi\left(\mathbb{E}[s]\right)$.
When designing basis functions, it is important to take the dimension $F$ \colr{of a basis function into account}. \colr{A larger dimension} $F$ will increase the computation time, as the algorithm requires computing the inverse of \colr{an} $F \times F$ matrix $A$. A high-dimensional parametrisation also increases the danger of overfitting. On the other hand, decreasing $F$ may result in a \colr{worse} approximation of the value function as the features become too general to accurately represent the state space.
In Table \ref{tab:basis_functions}, four possible basis functions are proposed. Basis functions $1$ and $4$ are similar to those used \colr{in} \citet{Hulshof2016}.

\begin{table}[h!]
\caption{Selected basis functions for least-squares policy iteration}
\label{tab:basis_functions}
\begin{tabular}{llll}
\hline
   & \textbf{Feature}    & \textbf{Basis function}      & \textbf{$F$}  \\
\hline
1. & \begin{tabular}[t]{@{}l@{}}Number of patients in queue $j \in \mathcal{J}$ \\ waiting for $w\in\mathcal{W}_j$ time periods\end{tabular}   & $s_{j, w}$, $\forall j \in \mathcal{J}, w\in \mathcal{W}_j$   & $\sum_{j\in\mathcal{J}}|\mathcal{W}_j|$ \\
2. & \begin{tabular}[t]{@{}l@{}}Number of patients in queue $j \in \mathcal{J}$ \\who are  early, on time, and late \\for treatment\end{tabular} & \begin{tabular}[t]{@{}ll@{}}$\sum_{w \leq u_j - 1} s_{j, w}$,\\ $ s_{j, u_j}$,&$\forall j \in \mathcal{J}$\\ $\sum_{w > u_j} s_{j, w}$,\end{tabular} & $3 \cdot |\mathcal{J}|$       \\
3. & \begin{tabular}[t]{@{}l@{}}Cost of patients in queue $j \in \mathcal{J}$ \end{tabular}     & $\sum_{w \geq u_j} c_{j, w} \cdot s_{j, w}$, $\forall j \in \mathcal{J}$    & $|\mathcal{J}|$ \\
4. & \begin{tabular}[t]{@{}l@{}}Total number of patients in queue $j \in \mathcal{J}$\end{tabular} & $\sum_{w\in\mathcal{W}_j} s_{j, w}$, $\forall j \in \mathcal{J}$  & $|\mathcal{J}|$  \\
\hline
\end{tabular}
\end{table}

The dimension of the first basis function may be too large as it is \colr{the sum} of the \colr{cardinalities of sets $\mathcal{W}_j$} of patient \colr{waiting} times. This also forces us to upper bound the \colr{waiting} times, while this is not strictly necessary for the other basis functions. However, an upper bound on \colr{waiting} times is necessary in any case for the policy-LP, which becomes slower to solve if the upper bound is increased. Since the number of \colr{queues $|\mathcal{J}|$} is limited, the other basis functions should be manageable. It is also possible to use a combination of multiple basis functions, for example, the cost of each queue \textit{and} the total number of patients in each queue.

A term of $1$ should be added to each basis function, increasing the dimension by one. If we do not do this, the all-zero state will always have value $V(\mathbf{0}) = \theta^{\top} \phi(\mathbf{0}) = 0$, which will not necessarily correspond with the exact value for that state. This term will improve the ability to approximate the exact value function. However, it does not change the resulting policy, as it only adds a constant term to the approximate value function and therefore does not affect which action results in the maximum value.

If the exact value function $V^*$ is known, we can use it to determine which basis functions should be used for LSPI. Given a basis function $\phi$, linear regression can be used to calculate the $\theta$ minimizing the mean squared error (MSE) between the exact \colr{value function $V^*$} and approximate \colr{value function $\overline{V}=\theta^T\phi$} for a training set of $n$ states $\{s_1, \ldots, s_n\} \subset \mathcal{S}$ \citep{Powell2019}:
\begin{equation}\label{eq:MSE}
\theta^* = \colr{\arg}\min_{\theta \in \mathbb{R}^F} \left( \sum_{i=1}^n \left( V^*(s_i) - \theta^{\top} \phi(s_i) \right)^2 \right).
\end{equation}

After using a training set of $n$ states to find $\theta\colr{^*}$ for a given basis function, the performance of that basis function can be measured by calculating the MSE on a test set of $m$ states:
\begin{equation*}
\text{MSE}(\phi, \{s_1, \ldots, s_m\}) = \sum_{i=1}^m \left( V^*(s_i) - \theta^{* \top} \phi(s_i) \right)^2.
\end{equation*}
The training and test set should be disjoint, as the test set measures how well the parametrisation found can generalise to unseen samples. Furthermore, $k$-fold cross-validation can be used to obtain an accurate estimation of the true MSE for the whole dataset \citep{Koutroumbas2008}.  
This process should be repeated for different basis functions. The \colr{basis function} resulting in the lowest MSE provides the closest approximation to the exact value function\colr{, and thus, should be used for LSPI}.

\subsection{Integer linear programming}\label{sec:solution_lp}

This section presents the integer linear program which can be used to solve the deterministic, finite-time-horizon variant of the MDP presented in Section \ref{sec:model}. The same sets, parameters, and variables are used; a summary can be found in Table \ref{tab:model_svp}.
Some changes must be made to the model as an ILP is deterministic and cannot handle an infinite time horizon. \colr{To remove all stochasticity from the problem, we let $q_{i,j}$ represent transition fractions instead of transition probabilities. In addition, we take for $\lambda_{j,t}$ the expected number of patients arriving in queue $j\in \mathcal{J}$ at time $t\in\mathcal{T}$.} \colr{Instead of} $\mathcal{T} = \mathbb{Z}^+$, we use a finite time horizon $\mathcal{T} = \{0, 1, \ldots, T\}$ for some $T \in \mathbb{Z}^+$.

\colr{Note that the ILP is similar to the policy-LP \eqref{lp:policy}, except that we take multiple time periods into account. To account for less precision in the later time periods, we use discount factor $\gamma$. The initial state is fixed to the current state for which the ILP is solved.}

\begin{maxi!}<b>
{a \in \mathcal{A}^T} {\sum_{t \in \mathcal{T}} \gamma^t C(s_t, a_t)}{}{ }
\addConstraint{s_{j, w, 0}} {= s_{j,w}^{current}, } {\quad\forall j \in \mathcal{J}, w\in \mathcal{W}_j}
\addConstraint{s_{j, 0, t}} {= \lambda_{j, t} + \sum_{i \in \mathcal{J}}\sum_{w\in\mathcal{W}_j}  q_{i,j} a_{i,w, t\colr{-1}}, } {\quad\forall j \in \mathcal{J}, t \in \mathcal{T},t>0}
\addConstraint{s_{j, w, t}}{= s_{j, w-1, t-1} - a_{j, w-1, t-1},} {\quad\forall j \in \mathcal{J}, 1\leq w \leq W_j, t \in \mathcal{T},t>0}
\addConstraint{s_{j, W_j,t}}{ = \sum_{w = W_j-1}^{W_j} s_{j, w,t} - a_{j, w,t},}{\quad \forall j \in \mathcal{J}, t \in \mathcal{T}}
\addConstraint{a_{j, w, t}}{ \leq s_{j, w, t}, }{\quad \forall j \in \mathcal{J}, w \in\mathcal{W}_j, t \in \mathcal{T}}
\addConstraint{\sum_{j \in \mathcal{J}}  \sum_{w\in\mathcal{W}_j} \zeta_{j, r}a_{j, w,t}}{\leq \eta_{r,t} ,}{\quad \forall r \in \mathcal{R}, t \in \mathcal{T}}
\addConstraint{a_{j, w,t}}{\in \mathbb{Z}^+, \label{constr:ILP_integer}}{\quad \forall j \in \mathcal{J}, w\in\mathcal{W}_j, t \in \mathcal{T}.}
\end{maxi!}

To make the problem faster to solve, integer constraints \eqref{constr:ILP_integer} could be relaxed to $a_{j, w, t} \geq 0$. In that case, the actions in the optimal solution must be rounded down to the nearest integer to ensure feasibility, as in the policy-LP. Again, this could result in loss of optimality, and should therefore not be done in the event that the ILP is actually used for planning. 

The solution to the ILP is a vector $a = (a_0, \ldots, a_T)$ of actions to perform over a time horizon of length $T$. However, actually performing these actions may not lead to good results in the real world for two reasons. The first is that, since the solution is fixed and transition fractions are used rather than probabilities, the solution is not \colr{optimal for real-life realizations}. The second reason is that the finite time horizon distorts the solution towards the end, as future consequences of those actions are not taken into account.
It is more realistic that a capacity planning department making use of this ILP would solve it every time a planning decision should be made, with $s_0$ as the (predicted) state for which a schedule must be made, and only implement the first action $a_0$. We refer to this as a rolling time horizon. $T$ should still be set to a longer period of time, such as one year, so that future consequences of the first action are taken into account. 

\subsection{Decision rules}\label{sec:solution_decisionrule}
In this section, we formulate four decision rules that can provide a feasible solution to the timeslot allocation problem \colr{at a certain time $t\in \mathcal{T}$}. The following decision rules are considered. \colr{Note that the index $t$ is omitted as we consider the state and resource capacity at a certain fixed time $t\in\mathcal{T}$.}
\begin{itemize}
\item \textbf{Highest Contribution:} Treat a patient that would result in the highest increase in contribution if they are treated now. This is expressed as the sum of the reward and the cost of the patient; by treating that patient, we gain the associated reward and avoid the associated cost:
\begin{equation*}
(j, w) = \arg \max \left\{\colr{c_{j, w}} + r_{j} :\; j \in \mathcal{J}, w \in \mathcal{W}_j, s_{j, w} \geq 1 \right\} .
\end{equation*}
Continue until all resource capacity has been used or the queues are empty.
\item \textbf{Highest Cost Queue:} Treat a patient from the queue with the highest total cost, i.e., 
\begin{equation*}
j = \arg \max \left\{ \colr{\sum_{w\in\mathcal{W}_j} c_{j, w}} \cdot s_{j, w}  :\; j \in \mathcal{J}\right\} .
\end{equation*}
The patient \colr{with the highest waiting time in the selected queue $j\in\mathcal{J}$} is treated. Recalculate the total queue costs after \colr{treating} the patient and continue until all resource capacity has been used or the queues are empty.
\item \textbf{Longest Queue:} Treat a patient from the longest queue, i.e., 
\begin{equation*}
j = \arg \max \left\{ \sum_{w\colr{\in \mathcal{W}_j}} s_{j, w} :\; j \in \mathcal{J} \right\} .
\end{equation*}
The patient \colr{with the longest waiting time} in the selected queue $j\in\mathcal{J}$ is treated. Recalculate the queue length after selecting the patient and continue until all resource capacity has been used or the queues are empty.
\item \textbf{Split Cost:} 
\colr{For this decision rule, we assume that each queue $j\in \mathcal{J}$ only requires one resource $r_j\in\mathcal{R}$. To determine the number of patients $a_j$ to treat from queue $j\in \mathcal{J}$, we distribute the available resource capacity of resource $r_j\in\mathcal{R}$ over all queues $i\in \mathcal{J}_{r_j}$ that require resource $r_j$, proportional to the cost of the queue.} 
\begin{equation*}
a_{j} = \min \left( \colr{\sum_{w\in\mathcal{W}_j}s_{j,w}}, \floor*{\eta_{r_j} \cdot \frac{\sum_{w \colr{\in\mathcal{W}_j}} c_{j, w} \cdot s_{j, w}}{\sum_{i \in \mathcal{J}_{r_j}}  \sum_{w\colr{\in\mathcal{W}_j}} c_{i, w} \cdot s_{i, w}}}\right) , \qquad \forall j \in \mathcal{J},
\end{equation*}
\colr{Note that we cannot treat more than $\sum_{w\in\mathcal{W}_j}s_{j,w}$ patients for queue $j\in \mathcal{J}$.} Within each queue, \colr{the patients with the highest waiting time are treated}. 
\end{itemize}

These rules are simple and easy to implement and understand. For hospitals, the use of decision rules is beneficial as they do not require complex optimization software and a solution can be found quickly. Also, it is much easier to explain to hospital staff why certain decisions should be made. The rule itself provides an explanation, for example, ``we should always treat our most urgent patients first''. This is in contrast with the other approaches, which are essentially black-boxes, as \colr{mathematical knowledge is required to explain} why the solution provided is optimal. However, we expect that this increase in simplicity and explainability comes with a loss of performance.

\subsection{Static allocation} \label{sec:solution_static}
Finally, we consider the method currently used by SMK for timeslot allocation within OD sessions. The hospital generally does not make adaptations to the number of FA, FU, and DA appointments that are planned within an OD session. Rather, a static allocation is used with a fixed number of appointments for each appointment type, as this can be done quickly and easily. \colr{Since manually adapting a surgeon's timeslots is a lot of work, it is only done when an intervention is necessary to prevent overly long or short waiting lists.}

For most surgeons, there are three different static allocations \colr{used} to plan an OD session, \colr{depending on whether} the surgeon is aided by zero, one, or two assistants. Schedules for sessions with more assistants allow for more OD appointment slots. However, the proportion of FA, FU, and DA appointments within one session remains roughly the same. The hospital determines these proportions by calculating the ratios between the appointment types for previous patients of the surgeon, and total OD time. Under the assumption that future patients and care pathways are similar to those in the past, this should ensure that waiting lists remain relatively balanced. We apply this method as well, in order to provide an accurate comparison between the new solution methods discussed in this paper and the current allocation method of SMK.

%% file: data_analysis.tex
\section{\colr{Test instances}}\label{sec:data_analysis}
\colr{In order to test how well the solution methods from Section \ref{sec:solution_methods} perform on the model constructed in Section \ref{sec:model}, we collect real-life data to create three test instances.} \colr{Based on this real-life data}, we generate two artificial \colr{instances}: one small test \colr{instance} (Section \ref{sec:data_small}), and one large test \colr{instance} (Section \ref{sec:data_large})\colr{, which are used to define the parameters of the solutions methods.} Section \ref{sec:data_smk} discusses  the full-scale instance based on actual SMK data.

The two test instances are required for a number of reasons. First, the small test \colr{instance} is required as we need to have an \colr{instance} which is solvable using exact value iteration (EVI). This \colr{instance} is used to compare how well a policy obtained through LSPI can approximate the optimal policy found through EVI. Therefore, the small test \colr{instance} should have a state space small enough to permit an exact solution. Secondly, LSPI requires some parameter tuning, and for this purpose, it is desirable to have a test \colr{instance} which is similar to real-life, but not quite as large so that it is faster to solve and experiments can be performed efficiently. Due to the strict size constraint, the small test \colr{instance} is too different from the real-life problem, so a second, large test \colr{instance} is defined. For both the test \colr{instances} and the complete SMK instance, each time period represents a planning period of two weeks.

\subsection{\colr{Small test instance}}\label{sec:data_small}
\begin{figure}[h!]
	\centering
	\includegraphics[width=0.4\textwidth]{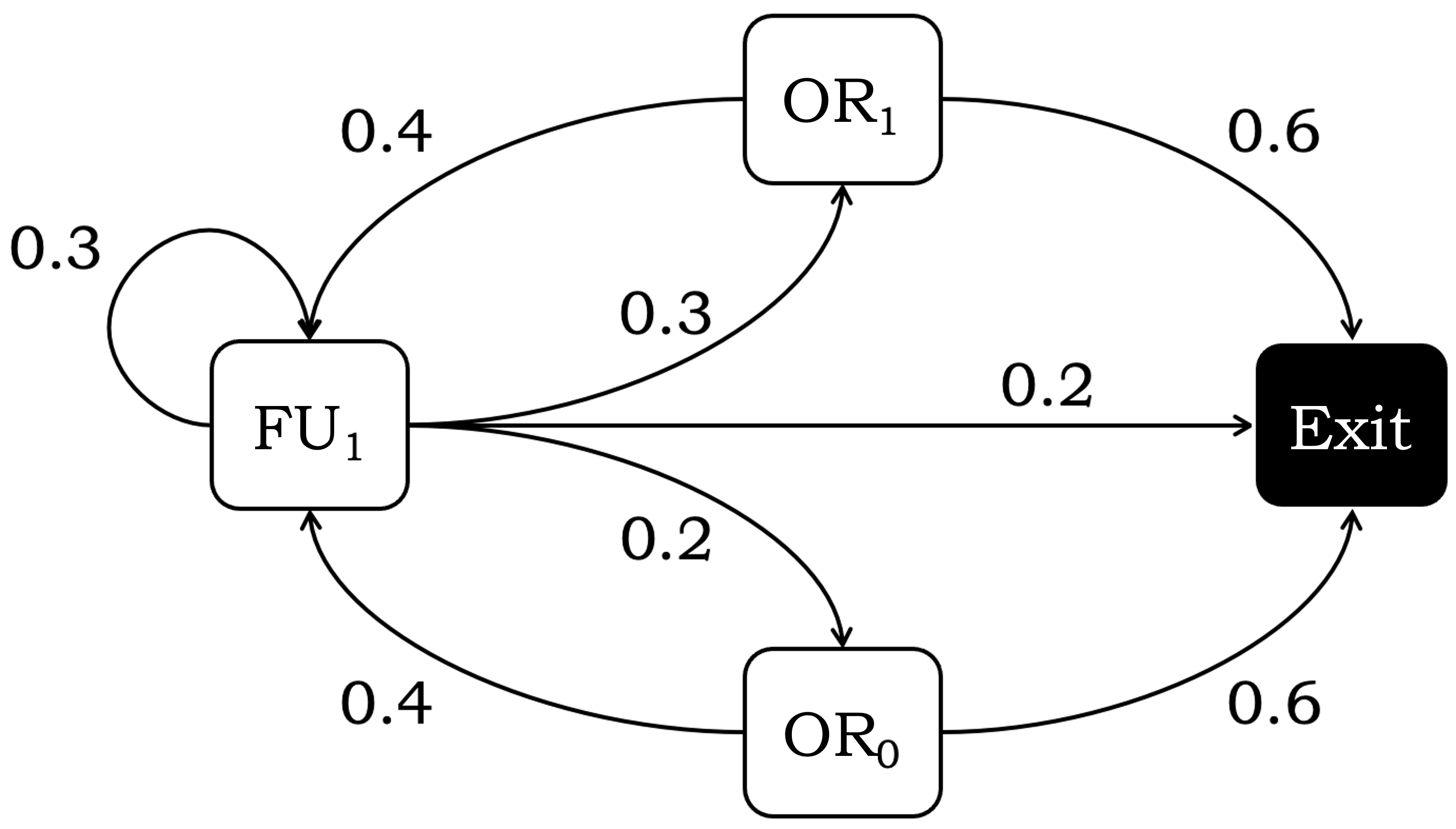}
	\caption{Queues and transition probabilities for the small test \colr{instance}}
	\label{fig:test_problem}
\end{figure}

In order to compare \colr{the policy
obtained through LSPI with the optimal policy found through EVI}, we generate a test \colr{instance} with a manageable state space. 
\colr{This small test instance has \colr{three} queues $\mathcal{J} = \{FU_1,\, OR_0,\, OR_1\}$, of which the FU queue has an access time target $u_{FU_1}$ of 1 and the OR queues have an access time target $u_{OR_0}$ and $u_{OR_1}$ of  0 and 1, respectively. The resource capacities of the resources $\mathcal{R}=\{OD,\, OR\}$ }are set to $\eta_{OD,t} = 4$ and $\eta_{OR,t} = 2$ for all time periods \colr{$t\in\mathcal{T}$.} The FU queue requires \colr{$\zeta_{FU_1,OD}=1$} unit of OD time per patient and the OR queue\colr{s} require \colr{$\zeta_{OR_0,OR}=\zeta_{OR_1,OR}=1$} unit of OR time per patient. \colr{The other resource requirements $\zeta_{FU_1,OR}$, $\zeta_{OR_0,OD}$, and $\zeta_{OR_1,OD}$ are 0.}  
The upper bound on \colr{the waiting} time, $W_j$, is \colr{set to} $2$ \colr{for each queue $j\in\mathcal{J}$}. 
The transition probabilities between queues are shown in Figure \ref{fig:test_problem}. Three new patients enter the FU queue in every time period, so $\lambda_{FU_1,t} = 3$ \colr{for all $t\in\mathcal{T}$} and \colr{$\lambda_{OR_0,t}$ and $\lambda_{OR_1,t}$ } are $0$ \colr{for all $t\in\mathcal{T}$}.
The cost for keeping a patient waiting \colr{in queue $j\in \mathcal{J}$ for $w\in \mathcal{W}_j$ time periods} is set to $c_{j, w} = \omega_j \cdot \frac{w}{u_j+1}$ where $\omega_{FU_1} = 1$ and $\omega_{OR_0}=\omega_{OR_1}  = 4$. Adding $1$ to $u_j$ in the denominator is necessary here as one of the \colr{access time targets} is $0$. The rewards for treating patients are $r_{FU_1} = 1$ and $\colr{r_{OR_0} =r_{OR_1}= 4}$.

\subsection{\colr{Large test instance}}\label{sec:data_large}
To be able to test different parameter settings and algorithm variations for the solution methods in an efficient manner, we create a test \colr{instance} that is a scaled-down but representative version of the full problem. The following parameters are used for the large test \colr{instance}. There are \colr{five} queues: $\mathcal{J} = \{FA_2,\, FU_4,\, OR_2,\, OR_4,\, DA_3\}$\colr{, for which the access time targets are given by $u_{FA_2}=2$, $u_{FU_4}=4$, $u_{OR_2}=2$, $u_{OR_4}=4$, and $u_{DA_3}=3$.} \colr{As for the small test instance, we have two resources given by $\mathcal{R}=\{OD,OR\}$.} The FA, FU and DA queues require one unit of OD time\colr{, i.e., $\zeta_{FA_2,OD}=\zeta_{FU_4,OD}=\zeta_{DA_3,OD}=1$} and the OR queue\colr{s} require one unit of OR time per patient\colr{, i.e., $\zeta_{OR_2,OR}=\zeta_{OR_4,OR}=1$}. \colr{The other resource requirements $\zeta_{FA_2,OR}$, $\zeta_{FU_4,OR}$, $\zeta_{OR_2,OD}$, $\zeta_{OR_4,OD}$, $\zeta_{DA_3,OR}$ are 0.} Resource capacities are set to $\eta_{OD,t} = 16$ and $\eta_{OR,t} = 2$ \colr{for all time periods $t\in\mathcal{T}$}. Transition probabilities are shown in Table \ref{tab:transitions_large_test}. In each row, the probability of a patient transferring from the appointment type in that row to the appointment type indicated in each column is shown. The `Exit' column is the probability that the patient will leave the system. 8 new patients enter the FA queue in every time period, so $\lambda_{FA_2,t} = 8$ \colr{for all $t\in\mathcal{T}$} and \colr{$\lambda_{FU_4,t}$, $\lambda_{OR_2,t}$, $\lambda_{OR_4,t}$, and $\lambda_{DA_3,t}$} are $0$ \colr{for all $t\in \mathcal{T}$}.

\begin{table}[h!]
\centering
\caption{Transition probabilities $\colr{q_{i,j}}$ for the large test \colr{instance}}
\label{tab:transitions_large_test}
\begin{tabular}{|c|c|c|c|c|c|c|}
\hline
       & $FA_2$ & $FU_4$ & $OR_2$ & $OR_4$ & $DA_3$ & Exit   \\ \hline
$FA_2$ & $0$    & $0.5$  & $0.01$  & $0.1$ & $0$    & $0.39$ \\ \hline
$FU_4$ & $0$    & $0.4$  & $0.02$  & $0.15$ & $0$  & $\colr{0.43}$ \\ \hline
$OR_2$ & $0$    & $0.2$  & $0$  & $0$ & $0.75$  & $0.05$ \\ \hline
$OR_4$ & $0$    & $0.25$  & $0$  & $0$    & $0.7$    & $0.05$  \\ \hline
$DA_3$ & $0$    & $0.6$  & $0$  & $0$    & $0$    & $0.4$  \\ \hline
\end{tabular}
\end{table}

\colr{We set the upper bound on access times to $W_{j} = 3 \cdot u_j$ for all $j \in \mathcal{J}$. This allows higher waiting times for queues in which patients generally wait longer due to higher access time targets.}

The cost for keeping a patient waiting \colr{in queue $j\in \mathcal{J}$ for $w\in \mathcal{W}_j$ time periods} is set to $c_{j, w} = \omega_j \cdot \frac{w}{u_j+1}$ where $\omega_{FA_2} = 2, \omega_{FU_4} = 1, \omega_{OR_2} = \omega_{OR_4} = 4$ and $\omega_{DA_3} = 1$. The rewards for treating patients are $r_{FA_2} = r_{FU_4} = 2, r_{OR_2} = r_{OR_4} = 10$ and $r_{DA_3} = 1$. 

\subsection{\colr{Case study instance}}\label{sec:data_smk}
In order to solve the timeslot allocation problem \colr{for a real-life instance}, we collected data at the SMK. This section describes how this data was collected, and which parameters are found from the data to specify \colr{the real-life} instance.

\colr{There are nine queues: $\mathcal{J} = \{FA_2,\, FU_3,\, FU_6,\, FU_{12},\, OR_1,\, OR_2,\, OR_4,\, OR_6,\, DA_3\}$, for which the access time targets are given by $u_{FA_2}=2$, $u_{FU_3}=3$, $u_{FU_6}=6$, $u_{FU_{12}}=12$, $u_{OR_1}=1$, $u_{OR_2}=2$, $u_{OR_4}=4$, $u_{OR_6}=6$, and $u_{DA_3}=3$. Again, we have two resources given by set $\mathcal{R}=\{OD,OR\}$.} The outpatient department resource capacity \colr{for the considered surgeon} is calculated by dividing the number of available OD timeslots in one year, by $26$ (the number of two-week periods in one year). This gives $\eta_{OD,t} = 121$ \colr{for all $t\in\mathcal{T}$}. Note that in real-life, the number of OD appointment slots \colr{varies over the year} as the number of OD dayparts scheduled varies per week. When using a solution method for real-time planning, the resource capacity can be adjusted per time period.
Similarly, the operating room resource capacity is calculated by dividing the number of surgeries performed by the considered surgeon in one year, by $26$. This gives $\eta_{OR,t} = 9$ \colr{for all $t\in\mathcal{T}$}.
First appointments take up two OD timeslots, while follow-up and discharge appointments only require one \colr{OD} timeslot. \colr{We assume} surgeries take up one OR timeslot. This gives $\zeta_{FA_2, OD} = 2$, and $\zeta_{FU_3, OD} =\zeta_{FU_6, OD} =\zeta_{FU_{12}, OD} = \zeta_{DA_3, OD} = \zeta_{OR_1, OR}=\zeta_{OR_2, OR}= \zeta_{OR_4, OR}= \zeta_{OR_6, OR} = 1$. \colr{The other resource requirements are 0.} \colr{As for the large test instance, the upper bounds on the waiting times are set to $W_{j} = 3 \cdot u_j$ for all $j\in\mathcal{J}$.}

Patient pathways are extracted from the hospital information system and used to generate vectors of appointment types, for example: $(FA_2, FU_3, FU_6, OR_2, DA_3)$. We use this data to estimate the transition probabilities. These are summarized in Table \ref{tab:transition_probs_full}. In each row, the probability of a patient transferring to the appointment type indicated in each column is provided. The `Start' row shows the probability of each appointment type being the patient's first appointment with the surgeon in question. The `Exit' column shows the probability of the treatment being finished.

\begin{table}[h!]
\caption{Transition probabilities \colr{$q_{i,j}$} for the \colr{case study instance}}
\centering
\label{tab:transition_probs_full}
\begin{tabular}{|l|l|l|l|l|l|l|l|l|l|l|}
\hline
& $FA_2$   & $FU_3$   & $FU_6$   & $FU_{12}$                     & $OR_1$   & $OR_2$   & $OR_4$   & $OR_6$   & $DA_3$   & Exit     \\ \hline
Start                                           & $0.7116$ & $0$      & $0.2138$ & $0$                           & $0.0185$ & $0.0026$ & $0.0049$ & $0.0264$ & $0.0220$ & $0$      \\ \hline
$FA_2$                                          & $0.0037$ & $0.2400$ & $0.1298$ & \multicolumn{1}{r|}{$0.1010$} & $0.0012$ & $0.0049$ & $0.0055$ & $0.0753$ & $0.0147$ & $0.4238$ \\ \hline
$FU_3$                                          & $0.0028$ & $0.1951$ & $0.1127$ & $0.0919$                      & $0.0038$ & $0.0104$ & $0.0189$ & $0.0994$ & $0.0170$ & $0.4479$ \\ \hline
$FU_6$                                          & $0.0085$ & $0.1575$ & $0.0840$ & $0.0953$                      & $0.0028$ & $0.0028$ & $0.0038$ & $0.0660$ & $0.0151$ & $0.5642$ \\ \hline
$FU_{12}$                                       & $0$      & $0.1535$ & $0.0930$ & $0.0605$                      & $0.0023$ & $0.0023$ & $0$      & $0.0628$ & $0.0070$ & $0.6186$ \\ \hline
$OR_1$                                          & $0$      & $0.2$    & $0.0333$ & $0.0167$                      & $0.0667$ & $0$      & $0$      & $0$      & $0.3833$ & $0.3$    \\ \hline
$OR_2$                                          & $0$      & $0.1$    & $0$      & $0$                           & $0.0333$ & $0$      & $0.0333$ & $0.0333$ & $0.6667$ & $0.1333$ \\ \hline
$OR_4$                                          & $0$      & $0.1522$ & $0.0435$ & $0.0435$                      & $0.0217$ & $0$      & $0$      & $0$      & $0.5870$ & $0.1522$ \\ \hline
$OR_6$                                          & $0$      & $0.1421$ & $0.0299$ & $0.0175$                      & $0.0050$ & $0$      & $0$      & $0.0099$ & $0.7182$ & $0.0773$ \\ \hline
$DA_3$                                          & $0.0021$ & $0.3080$ & $0.2089$ & $0.0654$                      & $0$      & $0.0021$ & $0.0021$ & $0.0232$ & $0.0105$ & $0.3776$ \\ \hline
\end{tabular}
\end{table}

The number of new patients \colr{is} calculated by dividing the budgeted number of new patients that the surgeon should \colr{see} in one year by $26$. This gives $\lambda_{FA_2,t} = 28$ \colr{for all $t\in\mathcal{T}$}. However, this is only an indication of the number of first appointments a surgeon should have throughout one year\colr{, while} Table \ref{tab:transition_probs_full} shows that it is not uncommon for a patient to start their treatment with an appointment that is not an FA. These are often patients that have been transferred from another surgeon. In our model, we still view these as new patients, since they originate from outside our one-surgeon system. We see that approximately $30\%$ of patients start with a non-FA appointment. Therefore, the total number of new patients should be $28 \cdot \frac{1}{0.7116} \approx 40$ \colr{leading to setting} $\lambda_{j,t} = 40 \cdot q_{Start, j}$\colr{ for all $j\in \mathcal{J}$ and $t\in \mathcal{T}$}.

The costs and rewards used in the contribution function should reflect the priorities of SMK and relative importance of different appointment types. SMK has indicated that it is especially important to meet the \colr{access time targets} for the OR queue. Using all available capacity is also more important for the OR queue, as OR time is much more expensive to waste than OD time. The hospital has also expressed that meeting \colr{access time targets} for follow-up appointments is more important than for first appointments. \colr{Given this,} costs and rewards are determined by repeatedly simulating the system with the rolling-horizon LP as solution method, and adapting costs and rewards until desired behaviour (as described above) is shown.  
Resulting from this trial-and-error method, we obtained \colr{for each queue $j\in\mathcal{J}$ and waiting time $w\in\mathcal{W}_j$} the cost function $c_{j, w} = \omega_j \cdot \frac{w}{u_j}$ for $w \geq u_j$ and $0$ elsewhere, with $\omega_{FA_2} = 0.5$, $\omega_{FU_3} = \omega_{FU_6}=\omega_{FU_{12}}=3$, $\omega_{OR_1} = \omega_{OR_2}=\omega_{OR_4}=\omega_{OR_6}=10$ and $\omega_{DA_3} = 3$. Rewards are set to $r_{FA_2} = 5$, $r_{FU_3} =r_{FU_6} =r_{FU_{12}} = 3$, $r_{OR_1} =r_{OR_2} =r_{OR_4} =r_{OR_6} = 50$ and $r_{DA} = 3$.

%% file: results.tex
This section shows the results of \colr{applying the different solution methods} on the three considered instances. These experiments are aimed at investigating how the performance of the solution methods from Section \ref{sec:solution_methods} varies as different parameters of the algorithms are changed, and which settings result in the best performance. We also compare the performance of the solution methods with each other.

\colr{In Section \ref{sec:evaluation}, we describe the simulation method used to evaluate the performance of the various methods when used in practice.} Section \ref{sec:results_small} is devoted to the small test \colr{instance}, which is used to determine how well the approximate value function found using LSPI can approximate the true value function found using EVI, and which basis functions lead to the best approximation. In Section \ref{sec:results_large}, we test different LSPI parameters, explore the effects of changing $\gamma$ on the rolling-horizon LP, and compare the decision rules introduced in Section \ref{sec:solution_decisionrule}, on the large test \colr{instance}. In these sections, we assume that we have perfect information about the state for which we wish to select an action, i.e., we do not plan ahead. In Section \ref{sec:results_smk}, we apply the solution methods to the \colr{case study instance where we do consider planning ahead. The results suggest that a combination of the rolling-horizon LP and static allocation might give even better results, which is also explored in Section \ref{sec:results_smk}}. 

Data and code used to program solution methods and generate results in this section can be found at: \url{https://github.com/yannavdv/MSc_thesis_Yanna}.

\colr{\subsection{Evaluation method}\label{sec:evaluation}
To compare the performance of the different solution methods, we use simulation in which the solution methods are applied on a number of trials for a number of subsequent time periods. This requires to create patient pathways, which are used to generate an initial state and to generate newly arriving patients. 
For the SMK instance, we randomly draw patient pathways out of a historical dataset of realised patient pathways. Since there is no historical data for the small and large test instance, we generate patient pathways by starting with a first appointment (FA) and then randomly drawing the next stage or exit based on the transition probabilities. 

We initialize the simulation by drawing pathways for the selected number of initial patients and sampling their current stage and waiting time. Patients are modeled as tuples $\left< v, i, w \right>$ where $v$ is the patient pathway, $i$ is a counter indicating the stage of treatment the patient is in (starting at 1), and $w$ is the current waiting time. For example, patient $\left<(FA_4, FU_6, FU_3, OR_4, DA_3), 4, 2\right>$ indicates a patient waiting for surgery with an access time target of 4 who has been waiting for two time periods after completing one first appointment and two follow-up appointments. 

The current waiting time for the initial state is sampled from the typical distribution of waiting times for each appointment type at any moment in time. We have fitted a number of probability distributions (beta, gamma, exponential, truncated normal) to the SMK data and the exponential distribution resulted in the best fit. Based on this, we estimate that the waiting times of patients are exponentially distributed with scale parameter $\beta = u_j$, where $u_j$ is the access time target of queue $j\in \mathcal{J}$. The current waiting time is drawn from this distribution and rounded to the nearest integer less than or equal to $W_{j}$. 

The patient pathways that we simulate are unknown to the solution methods; the pathways are only used to determine the next state after performing an action. When the patient is `treated', $i$ is incremented by $1$ and $w$ is set to $0$. If a patient has to wait for another time period, $w$ is incremented by $1$ when $w$ was less than $W_j$ or stays $W_j$ when $w$ was already equal to $W_j$. When $i > |v|$ the patient has completed their treatment and is removed from the system. After each time period, a fixed number of new patients arrive for which we sample a patient pathway. These patients are initialized at stage 1 with a waiting time of 0, i.e., $i=1$ and $w=0$.

The performance of the solution methods is measured by whether patients are seen within the access times targets and capacity is used efficiently, see Section \ref{sec:results_access_times}. To efficiently find the best settings, we use an estimation for performance by calculating the average contribution per time period over a number of trials. After applying an action to the current state, we add the corresponding contribution to the total contribution of that trial. Within each trial, the same settings are used for each solution method, meaning that the initial patients and newly arriving patients in each time period are the same.}

\subsection{Results small test \colr{instance}}\label{sec:results_small}

The small test \colr{instance} is devised such that it is possible to find the exact solution using exact value iteration. \colr{To guarantee a bounded state space to be able to apply EVI, the queue lengths are bounded from above by 7 and 2 for the FU queue and OR queues, respectively.} EVI was used to calculate the exact value function, using $\gamma = 0.9$ and tolerance parameter $\epsilon = 0.1$. After seven iterations, taking over twelve hours, the convergence criterion was met. The resulting value function is used throughout the remainder of this section.

Now that the exact value function $V^*$ for the small test \colr{instance} is known, it can be used to measure which basis functions result in the closest approximation of $V^*$, using the linear regression method. The exact solution provides us with a full dataset of $373,248$ state-value pairs. $10$-fold cross-validation was performed on this dataset, comparing each of the proposed basis functions. The $10$ subsets are selected randomly. Pairwise combinations of basis functions were also evaluated. For several of these combinations, no optimal $\theta^*$ was found. These combinations are left out of the results. The average MSEs for the (combinations of) basis functions are shown in Table \ref{tab:basis_functions_mse}.

\begin{table}[h!]
\caption{Mean squared errors of basis functions for 10-fold cross-validation on small test \colr{instance}}
\centering
\label{tab:basis_functions_mse}
\begin{tabular}{lccccc}
\hline
 &  &                 \\
\textbf{Basis function} & $1$&$2$&$3$&$4$&$3 + 4$ \\
\textbf{Average MSE} &  $2.3055$ &$2.3054$&$4.0474$ &$12.0450$ & $3.5708$   \\
  &                      &                 \\
\hline
\end{tabular}
\end{table}

It is clear that basis functions 1 and 2 result in the lowest MSE. It should be noted that those functions are very similar for the small test \colr{instance}. \colr{For queues $FA_1$ and $OR_1$ having an access time target of 1}, the number of early patients is \colr{equal to} $s_{j, 0}$, the number of on time patients is \colr{equal to} $s_{j, 1}$ and the number of late patients is \colr{equal to} $s_{j, 2}$. The only difference between the basis functions arises for \colr{queue $OR_0$} with \colr{an access time target of} $0$, as it is impossible to be early in that queue, and patients with \colr{a waiting} time \colr{of} $1$ and $2$ \colr{time periods} are late. This explains why their MSE's are almost exactly the same.
It is interesting to see that, while basis functions $3$ and $4$ are of equal dimension, basis function $3$ results in a much better approximation of the exact value function. The cost of a patient \colr{$c_{j,w}$ in queue $j\in\mathcal{J}$ waiting for $w\in\mathcal{W}_j$ time periods} provides information that is important to be able to approximate the value of the state, as it gives an indication of the number of patients and their combined waiting times. 
Combining functions \colr{3 and 4} provides better results, but does  not perform as good as basis functions 1 and 2.

\colr{We now compare the performance of the approximate policy found through LSPI, the exact policy found through EVI, and the rolling-horizon LP on the small test instance with varying numbers of initial patients. For the rolling-horizon LP, we relax the integrality constraint as this reduces the computation time by a factor of 100 approximately. As we are performing many experiments throughout this and the following sections, in which the program must be solved thousands of times, we choose to relax the integer constraint during the remainder of this paper. This does come with loss of optimality, and therefore, we expect results to improve compared to the results shown in this paper when the integer constraint is applied.

We perform $50$ trials per number of initial patients with 30 subsequent time periods each, and discount factor $\gamma = 0.9$ for all methods. For the rolling-horizon LP, we set the time horizon to $T = 30$. For LSPI, we use basis function $1$ and dataset size $M = 5,000$.}

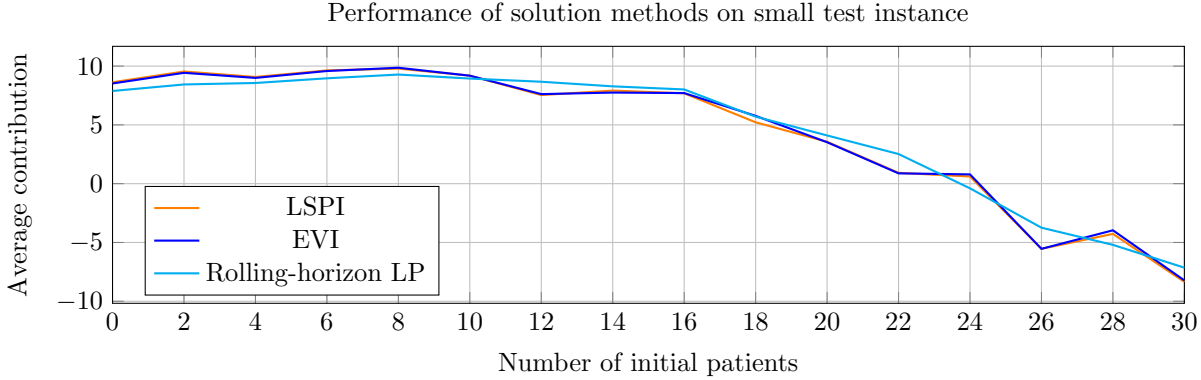
\begin{figure}[h!]
\begin{tikzpicture}
\begin{axis}[
	width=0.9\textwidth,
	height=5cm,
	title={Performance of solution methods on small test \colr{instance}},
	xlabel=Number of initial patients, 
	ylabel=Average contribution,
	xmin=0, xmax=30,
	legend pos=south west,
	xmajorgrids=true,
	ymajorgrids=true,
	every axis plot/.append style={thick}
]
\addplot[color=orange] coordinates {(0, 8.62533333333)
(2, 9.54066666667)
(4, 9.07466666667)
(6, 9.632)
(8, 9.78733333333)
(10, 9.19133333333)
(12, 7.538)
(14, 7.91666666667)
(16, 7.688)
(18, 5.22)
(20, 3.566)
(22, 0.92333333333)
(24, 0.608)
(26, -5.53866666667)
(28, -4.27)
(30, -8.35133333333)};
\addplot [color=blue] coordinates {(0, 8.52266666667)
(2, 9.42866666667)
(4, 8.996)
(6, 9.58266666667)
(8, 9.86066666667)
(10, 9.18466666667)
(12, 7.61733333333)
(14, 7.74933333333)
(16, 7.706)
(18, 5.766)
(20, 3.52466666667)
(22, 0.87933333333)
(24, 0.78933333333)
(26, -5.54066666667)
(28, -3.96466666667)
(30, -8.22666666667)};
\addplot[color=cyan]coordinates {(0, 7.88466666667)
(2, 8.43866666667)
(4, 8.56466666667)
(6, 8.96133333333)
(8, 9.282)
(10, 8.94133333333)
(12, 8.66933333333)
(14, 8.27866666667)
(16, 8.014)
(18, 5.70066666667)
(20, 4.102)
(22, 2.522)
(24, -0.39533333333)
(26, -3.74266666667)
(28, -5.204)
(30, -7.144)};
\legend{LSPI, EVI, Rolling-horizon LP}
\end{axis}
\end{tikzpicture}
\caption{\colr{Comparison of performance of LSPI, EVI, and rolling-horizon LP on the small test instance with varying numbers of initial patients.}} \label{fig:test_numpatients}
\end{figure}

Figure \ref{fig:test_numpatients} shows the results of this experiment. Recall that the  contribution should be maximized. We find that the performance of the policies resulting from LSPI and EVI is very similar. This implies that LSPI is able to approximate the exact value function well. \colr{Moreover}, the rolling-horizon LP produces quite similar \colr{and sometimes better} results \colr{than} LSPI and EVI. \colr{This can be explained by the fact that we bound the lenght of the queues when determining the optimal policy using EVI (resulting in an approximation) and that LSPI uses an approximation of the value function.} 
In the next section, we investigate \colr{various parameter settings to}  improve the LSPI algorithm. 

\subsection{Results large test \colr{instance}}\label{sec:results_large}
The large test \colr{instance} is constructed as described in Section \ref{sec:data_large}. The parameters reflect a scaled-down version of the parameters of the full \colr{instance} at SMK, \colr{thus enabling us} to perform experiments within reasonable time. The results of those experiments are described in this section. The static allocation method is \colr{not considered} here, as the optimal allocation for the large test \colr{instance} does not provide us with useful information about the optimal allocation for the full \colr{instance}. \colr{In Section \ref{sec:setLSPI}, we determine the values of $\gamma$, $\phi$, and $M$ which result in the best performance of LSPI. The value of $\gamma$ that results in the best performance for the rolling-horizion LP is determined in Section \ref{sec:setLP}. The best decision rule is determined in Section \ref{sec:decrules}.}

\subsubsection{Settings for LSPI}\label{sec:setLSPI}
For LSPI, only $\gamma$, $\phi$, and $M$ need to be chosen by the user (ignoring the constants $\epsilon$ and $\delta$, as these do not have a large influence on the algorithm as long as they are chosen small enough). Various methods have been proposed in literature, or can be devised by considering characteristics of the problem, to improve the solution quality and/or speed up convergence.

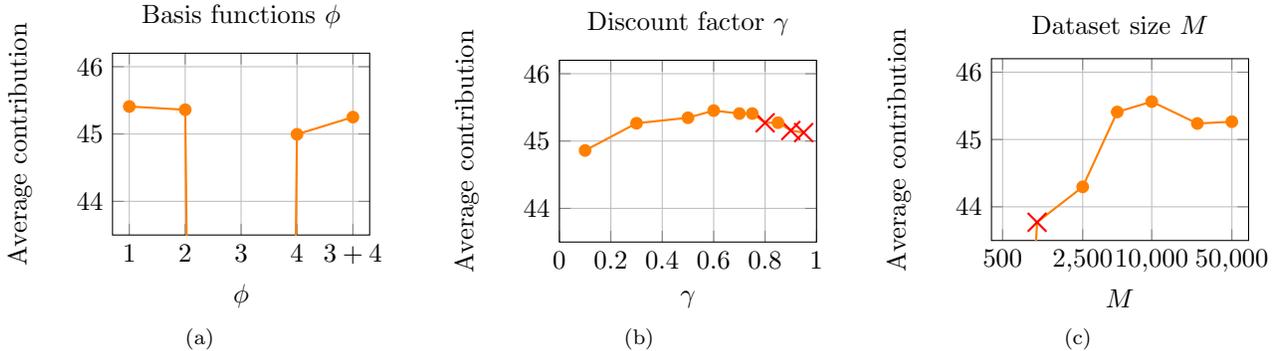
\begin{figure}[h!]
\subfloat[]{%
\begin{tikzpicture}
\begin{axis}[
	name=phi,
	width=5cm,
	height=4cm,
	log ticks with fixed point,
	xtick={1, 2, 3, 4, 5},
	xticklabels={$1$, $2$, $3$, $4$, $3+4$},
	title={Basis functions $\phi$},
	xlabel=$\phi$,
	ylabel=Average contribution,
	xmin=0.7, xmax=5.3,
	ymin=43.5, ymax=46.2,
	legend pos=south east,
	ymajorgrids=true,
	xmajorgrids=true,
	every axis plot/.append style={thick}
]
\addplot [color=orange]
coordinates {(1.0, 45.40857777777781)
(2.0, 45.36075555555555)
(3.0, -27.865899999999996)
(4.0, 44.99554444444445)
(5, 45.2513611111111)};
\addplot[red, mark=x, mark size=5pt, only marks] 
coordinates {(3, -27.865899999999996)};
\addplot[orange, mark=*, mark size=2pt, only marks] 
coordinates {(1, 45.40857777777781)
(2, 45.36075555555555)
(4, 44.99554444444445)
(5, 45.2513611111111)};
\end{axis}
\end{tikzpicture} \label{fig:phi}
} \hspace{0.5cm}
\subfloat[]{%
\begin{tikzpicture}
\begin{axis}[
	name=gamma,
	at={($(phi.east)+(3cm,0)$)},anchor=west,
	width=5cm,
	height=4cm,
	title={Discount factor $\gamma$},
	xlabel=$\gamma$,
	ylabel=Average contribution,
	xmin=0, xmax=1,
	ymin=43.5, ymax=46.2,
	legend pos=north west,
	ymajorgrids=true,
	xmajorgrids=true,
	every axis plot/.append style={thick}
]
\addplot[color=orange]
coordinates {(0.1, 44.86024444444445)
(0.3, 45.26452222222221)
(0.5, 45.345866666666666)
(0.6, 45.45241111111113)
(0.7, 45.40846666666669)
(0.75, 45.40857777777781)
(0.8, 45.271355555555566)
(0.85, 45.272600000000004)
(0.9, 45.156766666666655)
(0.95, 45.12593333333333)};
\addplot[red, mark=x, mark size=5pt, only marks] coordinates {(0.8, 45.271355555555566)
(0.9, 45.156766666666655)
(0.95, 45.12593333333333)};
\addplot[orange, mark=*, mark size=2pt, only marks] 
coordinates {(0.1, 44.86024444444445)
(0.3, 45.26452222222221)
(0.5, 45.345866666666666)
(0.6, 45.45241111111113)
(0.7, 45.40846666666669)
(0.75, 45.40857777777781)
(0.85, 45.272600000000004)
};
\end{axis}
\end{tikzpicture}%
}\hspace{0.5cm}
\subfloat[]{%
\begin{tikzpicture}
\begin{semilogxaxis}[
	name=datasize,at={($(phi.south)-(0, 2.3cm)$)},anchor=north,
	width=5cm,
	height=4cm,
	log ticks with fixed point,
	xtick={500, 2500, 10000, 50000},
	title={Dataset size $M$},
	xlabel=$M$,
	ylabel=Average contribution,
	xmin=400, xmax=70000,
	ymin=43.5, ymax=46.2,
	legend pos=south east,
	ymajorgrids=true,
	xmajorgrids=true,
	every axis plot/.append style={thick}
]
\addplot [color=orange]
coordinates {(500.0, 34.65407777777777)
(1000.0, 43.77044444444446)
(2500.0, 44.296666666666646)
(5000.0, 45.40857777777781)
(10000.0, 45.56277777777778)
(25000.0, 45.237777777777765)
(50000.0, 45.264)};
\addplot[red, mark=x, mark size=5pt, only marks] coordinates {(500.0, 34.65407777777777)
(1000.0, 43.77044444444446)};
\addplot[orange, mark=*, mark size=2pt, only marks] 
coordinates {(2500.0, 44.296666666666646)
(5000.0, 45.40857777777781)
(10000.0, 45.56277777777778)
(25000.0, 45.237777777777765)
(50000.0, 45.264)};
\end{semilogxaxis}
\end{tikzpicture}%
}
\caption{Performance of LSPI on large test \colr{instance} for different parameter settings. Settings for which there was no convergence are marked with a red cross.} \label{fig:test_lspi}
\end{figure}

\colr{We evaluate the performance of different settings for LSPI on the large test instance with varying numbers of initial patients. We perform $50$ trials per number of initial patients (drawn uniformly from the range $[50, 70]$) with the time horizon set to $30$ and, unless stated otherwise, discount factor $\gamma = 0.75$, basis function 1, and a dataset size $M = 5,000$.}
The same dataset is used for all experiments, except for the experiments on the size of the dataset. Since experiments have shown that the algorithm usually converges within approximately $4$ to $6$ iterations on the large test \colr{instance} for similar parameter settings, an upper bound is set on the number of iterations of $N = 20$. If the convergence criterion has not been met by this point, we assume that the algorithm will never converge for the settings used. \colr{Figure \ref{fig:test_lspi} shows the average contribution of the LSPI for different settings}. 

Figure \ref{fig:test_lspi}(a) shows the performance of the four different basis functions. LSPI using basis function $3$ (the cost of patients in each queue) does not converge.
Despite this, the combination of basis functions $3$ and $4$ does converge, and results in higher average contribution than basis function $4$ alone. LSPI converged for all other basis functions in 4-6 iterations, with basis function $1$ resulting in the best performance, and $4$ in the worst. For this reason, basis function $1$ (essentially the full state vector) is used for LSPI on the case study instance.

Figure \ref{fig:test_lspi}(b) shows the performance for different discount factors. The discount factor $\gamma$ determines to what extent we take into account future costs and rewards. If $\gamma = 0$, the policy is only based on the contribution gained by the current state and action, while for $\gamma \rightarrow 1$ we try to look infinitely far into the future. Intuitively, this tells us that increasing $\gamma$ will improve the policy, but also makes the value function much more difficult to approximate and can result in no convergence. This intuition is confirmed by the experiments testing $\gamma$, where we see that the \colr{average} contribution increases as $\gamma$ increases, until LSPI stops converging for $\gamma \geq 0.8$. Choosing $\gamma \approx 0.6$ results in the best policies.

Figure \ref{fig:test_lspi}(c) shows the performance for different dataset sizes. Increasing the size of the dataset that $\theta$ is trained on, $M$, should result in better policies as the algorithm can explore a larger part of the state space. The trade-off here is that this will increase the running time of the algorithm (linearly). The results show that LSPI will converge reliably for $M \geq 2{,}500$; the algorithm did not converge for $M = 500$ and $M = 1{,}000$. The average contribution increases as $M$ increases until $M = 10{,}000$. \colr{For values of $M$ higher than $10{,}000$}, the average contribution is lower. This could be a sign of overfitting, or it could indicate that, above a certain dataset size, the performance of LSPI is more sensitive to the contents of the dataset than to the size. 

\subsubsection{Setting $\gamma$ for the rolling-horizon LP}\label{sec:setLP}

\begin{figure}[h!]
\begin{tikzpicture}
\begin{axis}[
	name=lp_gamma,
	width=7.5cm,
	height=4cm,
	title={Performance of rolling-horizon LP versus $\gamma$},
	title style={yshift=0.3cm},
	xlabel=$\gamma$,
	ylabel=Average contribution,
	xmin=0, xmax=1,
	ymajorgrids=true,
	xmajorgrids=true,
	every axis plot/.append style={thick}
]
\addplot[color=cyan, x=gamma, y=contribution]
 coordinates {
(0.0, 43.66017222222221)
(0.05, 44.38270555555556)
(0.1, 44.21900555555555)
(0.15, 44.28484444444445)
(0.2, 44.12456666666666)
(0.25, 44.44562222222221)
(0.3, 44.594466666666676)
(0.35, 44.428888888888885)
(0.4, 44.55588333333335)
(0.45, 44.48743333333332)
(0.5, 44.525483333333334)
(0.55, 44.66024444444446)
(0.6, 44.466177777777794)
(0.65, 44.50538333333334)
(0.7, 44.54319444444443)
(0.75, 44.521483333333364)
(0.8, 44.51337777777778)
(0.85, 44.61821666666665)
(0.9, 44.62535000000002)
(0.95, 44.55435555555557)
(1.0, 44.029761111111114)};

\end{axis}

\begin{axis}[
	at={($(lp_gamma.east)+(2.3cm,0)$)},anchor=west,
	width=7.5cm,
	height=4cm,
	title={\colr{Computing time} rolling-horizon LP versus $\gamma$},
	title style={yshift=0.3cm},
	xlabel=$\gamma$,
	ylabel=\colr{Computing time (seconds)},
	xmin=0, xmax=1,
	xmajorgrids=true,
	ymajorgrids=true,
	every axis plot/.append style={thick}
]
\addplot[color=cyan]
[x=gamma, y=time] coordinates {(0.0, 0.0276905807)
(0.05, 0.02968495518)
(0.1, 0.03037063067)
(0.15, 0.03052774184)
(0.2, 0.0306142638)
(0.25, 0.03044340385)
(0.3, 0.03022244977)
(0.35, 0.03019725091)
(0.4, 0.03006590706)
(0.45, 0.03007730027)
(0.5, 0.03036206136)
(0.55, 0.03056937415)
(0.6, 0.03040029736)
(0.65, 0.0301037165)
(0.7, 0.03007007517)
(0.75, 0.03000246919)
(0.8, 0.03011733947)
(0.85, 0.03042427751)
(0.9, 0.03062827253)
(0.95, 0.03041713053)
(1.0, 0.03026515824)
};
\end{axis}
\end{tikzpicture}
\caption{Performance and \colr{computing} time of the rolling-horizon LP for different values of $\gamma$} \label{fig:test_lp_gamma}
\end{figure}
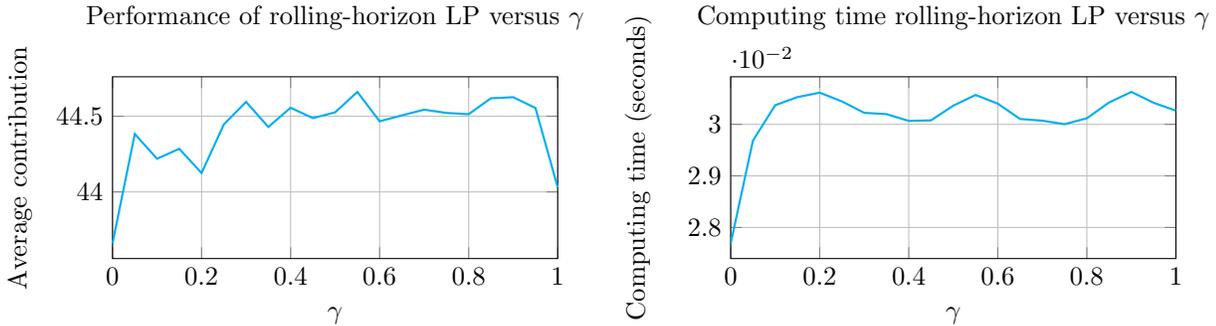

The discount factor $\gamma$ is, next to LSPI, also part of the objective function of the rolling-horizon LP. 
We evaluate the average contribution per time period and \colr{computing} time of the rolling-horizon LP for the large test \colr{instance}, where the average is taken over $50$ trials \colr{for different values of $\gamma$}. The time horizon\colr{, number of time periods,} and number of initial patients are all set to $30$. Results are shown in Figure \ref{fig:test_lp_gamma}, which shows that initially, \colr{the average} contribution and \colr{computing time} increase as $\gamma$ increases. Although there are minor fluctuations, the region $\gamma \in [0.3, 0.95]$ \colr{shows a stable performance, while} there is a large drop in performance from $\gamma = 0.95$ to $1$. The discount factor does not appear to have much influence on \colr{computing} time for $\gamma \geq 0.1$.

\subsubsection{Comparison of decision rules}\label{sec:decrules}

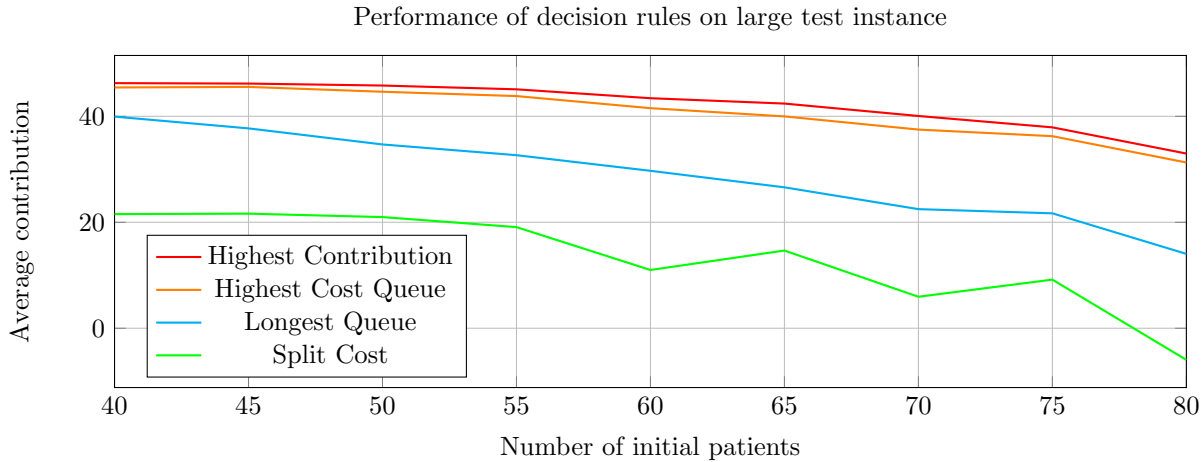
\begin{figure}[h!]
\begin{tikzpicture}
\begin{axis}[
	width=0.9\textwidth,
	height=6cm,
	title={Performance of decision rules on large test \colr{instance}},
	xlabel=Number of initial patients,
	ylabel=Average contribution,
	xmin=40, xmax=80,
	legend pos=south west,
	xmajorgrids=true,
	ymajorgrids=true,
	every axis plot/.append style={thick}
]
\addplot[color=red]
coordinates {(40.0, 46.25982222222221)
(45.0, 46.17706666666666)
(50.0, 45.80677777777778)
(55.0, 45.093355555555554)
(60.0, 43.41282222222222)
(65.0, 42.40397777777777)
(70.0, 40.06051111111112)
(75.0, 37.914644444444434)
(80.0, 32.97186666666667)
};
\addplot[color=orange]
coordinates {(40.0, 45.44393333333334)
(45.0, 45.548088888888884)
(50.0, 44.63655555555556)
(55.0, 43.81873333333333)
(60.0, 41.536466666666676)
(65.0, 39.98086666666666)
(70.0, 37.49666666666667)
(75.0, 36.25022222222222)
(80.0, 31.283155555555563)};
\addplot[color=cyan]
coordinates {(40.0, 39.94177777777778)
(45.0, 37.72057777777777)
(50.0, 34.687866666666665)
(55.0, 32.662977777777776)
(60.0, 29.71288888888889)
(65.0, 26.59506666666667)
(70.0, 22.480177777777776)
(75.0, 21.689111111111114)
(80.0, 14.034977777777781)};
\addplot[color=green]
coordinates {(40.0, 21.539466666666666)
(45.0, 21.6278)
(50.0, 20.980955555555553)
(55.0, 19.086000000000002)
(60.0, 10.997999999999998)
(65.0, 14.663711111111109)
(70.0, 5.9452)
(75.0, 9.173644444444443)
(80.0, -5.969466666666665)
};
\legend{Highest Contribution, Highest Cost Queue, Longest Queue, Split Cost}
\end{axis}
\end{tikzpicture}
\caption{Comparison of the decision rules for different numbers of initial patients.} \label{fig:test_decisionrule_large}
\end{figure}

The four decision rules described in Section \ref{sec:solution_decisionrule} were implemented and tested on the large test instance. Figure \ref{fig:test_decisionrule_large} shows the average contribution per time period gained over $50$ trials when executing each of the decision rules for 30 subsequent time periods. The results show that the Highest Contribution rule performs best, followed closely by the Highest Cost Queue rule. The Split Cost rule performs worst. 

The Highest Contribution rule is a direct, greedy maximization of the contribution function within each time period, without taking future effects of the decisions into account. Technically, the rule is equivalent to the rolling-horizon LP with $\gamma = 0$, as in that case the action is selected which maximizes the contribution in only the current time period. Because of this, we expect the performance of the decision rule to be slightly worse than that of the rolling-horizon LP.

%% file: results_smk.tex
\subsection{Results case study}\label{sec:results_smk}
Now that we have determined which hyperparameters to use for LSPI and the rolling-horizon LP, and which decision rule is best, we can apply the solution methods to the full timeslot allocation problem at SMK \colr{as described in Section \ref{sec:data_smk}}. 
\colr{In Section \ref{sec:comparison_full}, we compare the performance of the solution methods from Section \ref{sec:solution_methods} with and without planning ahead.}
\colr{Based on these results, we propose a hybrid method that combines the simplicity of the static allocation and the performance of the rolling-horizon LP in Section \ref{sec:hybrid}. Finally, in Section \ref{sec:results_access_times}, we compare how well the LP, static allocation, and hybrid method adhere to the access time targets and what the resulting unused capacity is}. 

\subsubsection{Comparison of \colr{solution methods}}\label{sec:comparison_full}
We compare the performance of the LSPI, the rolling-horizon LP, the Highest Contribution decision rule, and the static allocation used by SMK. We compare the average contribution over $100$ trials, \colr{consisting of} $30$ subsequent time periods each, and the following settings for the solution methods:

\begin{itemize}
\item \textbf{LSPI:} Basis function $1$, $\gamma = 0.6$, and dataset size $M = 50,000$. Due to the increased size of the problem, we expect a larger dataset than for the large test \colr{instance} is necessary to achieve good results. The LSPI algorithm took $3$ iterations, lasting $4$ hours, to converge to a parameter vector.
\item \textbf{Rolling-horizon LP:} Time horizon $T = 26$ (one year) and $\gamma = 0.75$ as we found that this discount factor performs better when planning ahead.
\item \textbf{Decision rule:} Highest Contribution.
\item \textbf{Static allocation:} The data suggests that approximately $50\%$ of OD time should go to first appointments, $42.5\%$ to follow-up appointments and $7.5\%$ to discharge appointments. \colr{For each time period $t\in\mathcal{T}$, this translates to treating 30 patients from the $FA_2$ queue, treating 17 patients from queues $FU_3$, $FU_6$, and $FU_{12}$ each, and treating 9 patients from queue $DA_3$} (recall that first appointments take twice as long as follow-up and discharge appointments). All OR capacity should be used, if possible. Within each queue, appointments are allocated to patients with the highest cost.
\end{itemize}

\begin{figure}[h!]
\begin{tikzpicture}
\begin{axis}[
	width=0.9\textwidth,
	height=6cm,
	title={Performance of solution methods on \colr{the case study instance}},
	xlabel=Number of initial patients,
	ylabel=Average contribution, 
	xmin=400, xmax=900,
	legend pos=south west,
	ymajorgrids=true,
	xmajorgrids=true,
	every axis plot/.append style={thick}
]
\addplot[color=orange]
coordinates {(400, 776.7389743589745)
(450, 773.2640064102565)
(500, 769.4118589743592)
(550, 756.8937820512821)
(600, 738.4841346153848)
(650, 711.8097115384614)
(700, 671.951057692308)
(750, 631.1358333333336)
(800, 582.3762500000005)
(850, 534.7434935897435)
(900, 483.81403846153853)};
\addplot[color=cyan]
coordinates {(400, 774.1401282051285)
(450, 774.556538461538)
(500, 774.7217628205127)
(550, 767.874455128205)
(600, 758.4765384615388)
(650, 745.1826602564101)
(700, 718.0250961538463)
(750, 686.2104487179485)
(800, 644.1457051282049)
(850, 601.4726282051279)
(900, 551.264230769231)};
\addplot[color=red]
coordinates {(400, 756.6809935897437)
(450, 748.8875000000002)
(500, 744.8829807692308)
(550, 739.95782051282)
(600, 735.335608974359)
(650, 727.994487179487)
(700, 710.9396153846153)
(750, 685.2771153846153)
(800, 647.1381730769232)
(850, 604.8461217948719)
(900, 556.6556089743589)};
\addplot[color=green]
coordinates {(400, 710.3047435897435)
(450, 686.8433653846156)
(500, 666.6329166666666)
(550, 639.1860897435897)
(600, 613.723044871795)
(650, 588.5219230769231)
(700, 559.1435256410255)
(750, 531.0309294871795)
(800, 506.7879166666669)
(850, 480.00554487179494)
(900, 445.9332051282054)};
\legend{LSPI, Rolling-horizon LP, Highest Contribution, Static allocation}
\end{axis}
\end{tikzpicture}
\caption{Comparison of policies resulting from LSPI, the rolling-horizon LP, the Highest Contribution decision rule, and static allocation on the case study instance when the number of initial patients is varied.} \label{fig:test_numpatients_full}
\end{figure}
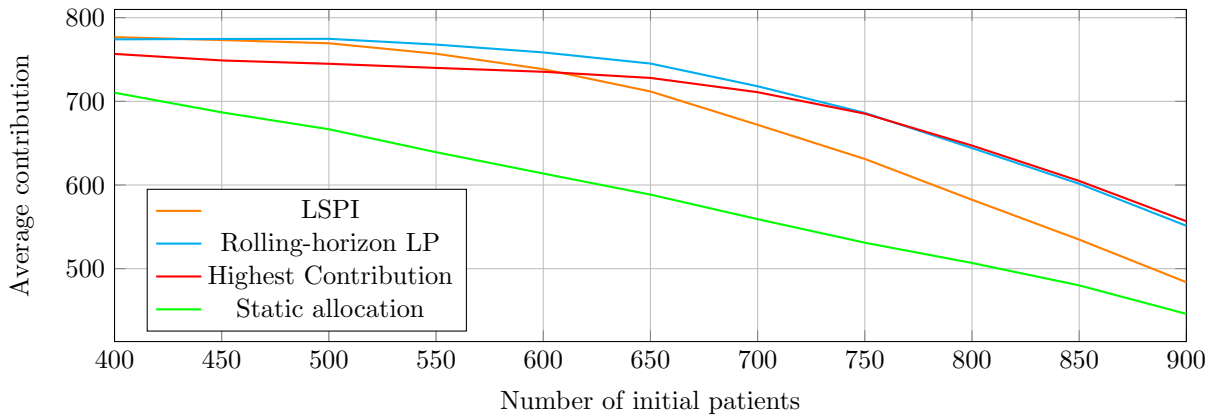

\colr{First, we assume that we have perfect information about the state for which we must choose an action, i.e., we do not plan ahead. Figure \ref{fig:test_numpatients_full} shows the results for varing numbers of initial patients.} 
The results show that for $400$ initial patients, LSPI results in the highest contribution, although only marginally higher than the rolling-horizon LP. The performance of LSPI quickly drops below the performance of rolling-horizon LP and decision rule as the number of initial patients is increased. The rolling-horizon LP steadily achieves good results, although it is slightly overtaken by the decision rule for over $750$ patients. The lowest contribution is accumulated by the static allocation.
The fact that the Highest Contribution decision rule performs slightly better than the rolling-horizon LP for a high number of patients is interesting, as the decision rule is essentially the same as the LP when we set $\gamma = 0$. This would suggest that with a higher number of patients, it is better to decrease $\gamma$ in order to deal with the high waiting lists immediately. This could be exploited through a scheme where $\gamma$ is adapted depending on the expected number of patients.
 
\colr{Second, we} investigate what happens in the more realistic situation that timeslots must be allocated to appointment types several weeks in advance, and a prediction of the future state must be used to make this scheduling decision, as explained in Section \ref{sec:model_planning_ahead}. Since our prediction of the future state will differ from the actual state, we expect to lose some performance quality. 

\begin{figure}[h!]
\begin{tikzpicture}
\begin{axis}[
	width=0.9\textwidth,
	height=6.5cm,
	title={Performance of solution methods when planning ahead},
	xlabel= Number of time periods to plan ahead ($p$),
	ylabel=Average contribution, 
	xmin=0, xmax=8,
	legend pos=south west,
	ymajorgrids=true,
	xmajorgrids=true,
	every axis plot/.append style={thick}
]
\addplot[color=orange]
coordinates {(0, 634.5825641025637)
(1, 579.9091666666669)
(2, 575.7975641025641)
(3, 572.0023717948718)
(4, 565.3469871794871)
(5, 551.2378205128203)
(6, 532.7303525641024)
(7, 514.5277243589742)
(8, 499.051346153846)};
\addplot[color=cyan]
coordinates {(0, 680.4366987179486)
(1, 682.1066346153847)
(2, 680.9666666666667)
(3, 679.0859294871797)
(4, 662.4032051282052)
(5, 657.5044551282051)
(6, 655.0693269230769)
(7, 646.9649038461534)
(8, 642.111730769231)};
\addplot[color=red]
coordinates {(0, 670.0786217948719)
(1, 670.4131089743589)
(2, 668.4742307692308)
(3, 666.3222435897435)
(4, 649.2467628205129)
(5, 643.0460576923077)
(6, 638.0756410256408)
(7, 627.6569871794873)
(8, 620.5774358974355)
};
\addplot[color=green]
coordinates {(0, 548.4534615384612)
(1, 548.4534615384612)
(2, 548.4534615384612)
(3, 548.4534615384612)
(4, 548.4534615384612)
(5, 548.4534615384612)
(6, 548.4534615384612)
(7, 548.4534615384612)
(8, 548.4534615384612)};
\legend{LSPI, Rolling-horizon LP, Highest Contribution, Static allocation}\end{axis}
\end{tikzpicture}
\caption{Comparison of policies resulting from LSPI, the rolling-horizon LP, the Highest Contribution decision rule and static allocations on the \colr{case study instance} when planning $p$ time periods ahead.} \label{fig:test_12weeks_full}
\end{figure}
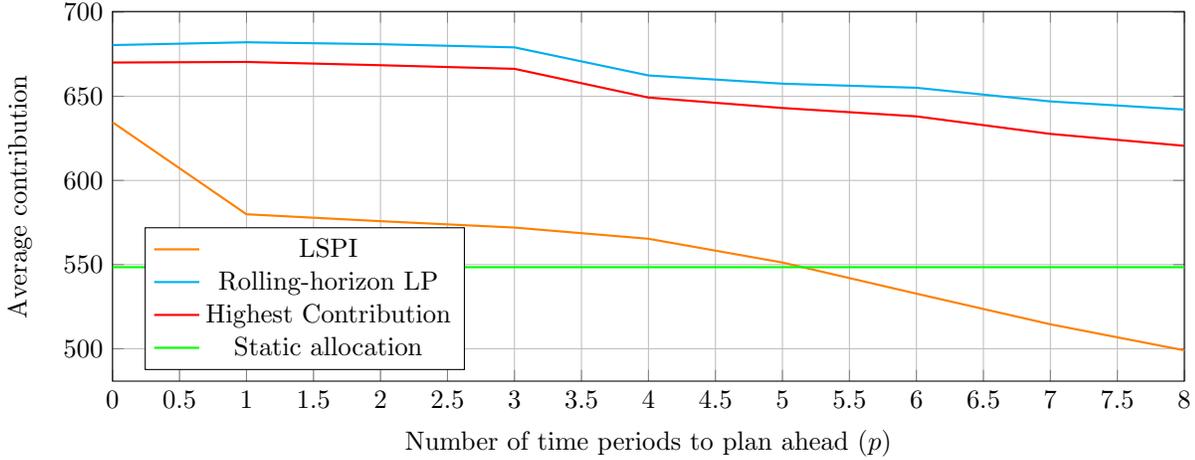

Figure \ref{fig:test_12weeks_full} \colr{compares the average} contribution of the different solution methods for different values of $p$. 
The number of initial patients in the simulation is drawn in each trial from $\mathcal{N}(700,200)$.  
When planning $p$ weeks ahead, in each time period, we calculate the predicted state for the current time period, $\widetilde{s}_t$, using the true state $p$ time periods ago, $s_{t-p}$, and the actions performed in the meantime. The first $p + 1$ time periods are ignored in the average contribution, as this is a warmup-period where $s_{t-p}$ does not yet exist. 

We see that, as expected, performance of the rolling-horizon LP, LSPI, and the decision rule decrease as the planning horizon increases, while that of the static allocation remains the same. \colr{The performance of LSPI decreases sharply} from $p = 0$ to $p = 1$ and \colr{then gradually worsens until it performs worse} than the static allocation method. It is interesting that the rolling-horizon LP works much better, while it is designed to solve the same problem as LSPI. This suggests that the linear parametrisation used in LSPI is too much of a simplification of the real-life problem. The decision rule performs quite well consistently, although it is not as good as the LP. The performance of the LP and Highest Contribution decision rule does not decrease much from $p=0$ to $p=3$.

\colr{The results show that the performance of the dynamic solution methods decreases as the planning horizon increases. SMK currently plans appointments and surgeries up to 12 weeks ($p=6$) in advance, which is appreciated by patients. Therefore, we investigate whether it is possible to statically allocate a percentage of timeslots 12 weeks in advance and allocate the remainder dynamically at a later point in time.}

\subsubsection{Hybrid method} \label{sec:hybrid}
\colr{We investigate a} hybrid method \colr{that} combines the static allocation method that plans $12$ weeks ahead, with the rolling-horizon LP that plans less far ahead. 
This allows us to make some appointments far enough in advance, while exploiting the improved knowledge we have of the predicted state when planning less far ahead. The static allocation method is used to determine the allocation of $\alpha\%$ of the timeslots $p=6$ time periods ($12$ weeks) before the time period to be scheduled. This is the \textit{fixed} allocation portion. Then, we use the rolling-horizon LP $\tau<p$ time periods prior to the week to be scheduled, to determine the allocation of the \colr{remaining $(1-\alpha)\%$ of the} timeslots, making up the \textit{dynamic} schedule portion.

Let $\overline{a}_j$ be the number of patients to treat from queue $j\in\mathcal{J}$, as determined by the static allocation. Let $K$ be some large number \colr{representing an upper bound on the queue length over all queues $j\in \mathcal{J}$}. Recall that $s_{j, w, t}$ denotes the sub-state size and $a_{j, w, t}$ the number of patients to treat from that sub-state. \colr{Constraints} \eqref{firsthybrid}-\eqref{lasthybrid} should be added to the rolling-horizon LP to \colr{ensure that the correct action is taken}.
\colr{Constraints} \eqref{firsthybrid} and \eqref{secondhybrid} determine \colr{the minimum of the number of patients $\sum_{w\in\mathcal{W}_j} s_{j, w, t}$ in queue $j\in\mathcal{J}$ at time period $t\in\mathcal{T}$ and the number of patients $\overline{a}_j \cdot \frac{\alpha}{100}$ to treat according to the fixed allocation portion}. \colr{For this, binary variables $y_{j,t}$ are introduced for each queue $j\in\mathcal{J}$ and time period $t\in\mathcal{T}$, which are 1 when the minimum is equal to $\sum_{w\in\mathcal{W}_j} s_{j, w, t}$ and 0 when the minimum is equal to $\overline{a}_j \cdot \frac{\alpha}{100}$.} \colr{Then, Constraints} \eqref{thirdhybrid} and \eqref{fourthhybrid} ensure that the action $\sum_{w\in\mathcal{W}_j} a_{j, w, t}$ chosen for \colr{queue $j\in\mathcal{J}$ at time period $t\in\mathcal{T}$}, is \colr{greater than or equal to this} minimum. This is necessary because simply forcing the action to be larger than \colr{or equal to} the static part of the allocation can lead to infeasibility, if the size of the queue is smaller in that time period.

\begin{align}
\overline{a}_j \cdot \frac{\alpha}{100} -  \sum_{w\in\mathcal{W}_j} s_{j, w, t} &\leq K \cdot y_{j, t} , \quad &\forall j \in \mathcal{J}, t \in \mathcal{T} , \label{firsthybrid}\\
\sum_{w\in\mathcal{W}_j} s_{j, w, t} - \overline{a}_j \cdot \frac{\alpha}{100} &\leq K \cdot (1-y_{j, t}) , \quad &\forall j \in \mathcal{J}, t \in \mathcal{T} , \label{secondhybrid}\\
\sum_{w\in\mathcal{W}_j} a_{j, w, t} &\geq \sum_{w\in\mathcal{W}_j} s_{j, w, t} - K \cdot (1-y_{j, t}) , \quad &\forall j \in \mathcal{J}, t \in \mathcal{T} \label{thirdhybrid}\\
\sum_{w\in\mathcal{W}_j } a_{j, w, t} &\geq \overline{a}_j \cdot \frac{\alpha}{100} - K \cdot y_{j, t} \quad &\forall j \in \mathcal{J}, t \in \mathcal{T} ,\label{fourthhybrid}\\
y_{j, t} &\in \{0, 1\} , \quad &\forall j \in \mathcal{J}, t \in \mathcal{T} .\label{lasthybrid}
\end{align}

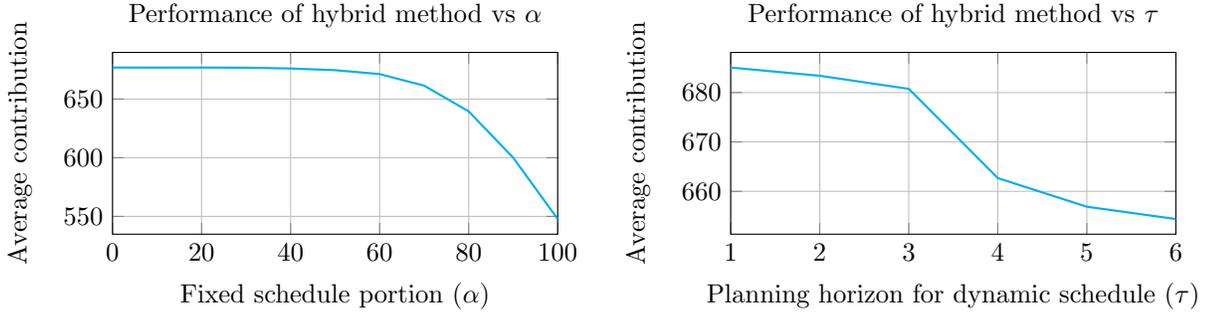
\begin{figure}[h!]
\begin{tikzpicture}
\begin{axis}[
	name=hybrid_alpha,
	width=7.5cm,
	height=4cm,
	title={Performance of hybrid method vs $\alpha$},
	xlabel=Fixed schedule portion ($\alpha$),
	ylabel=Average contribution,
	xmin=0, xmax=100,
	xmajorgrids=true,
	ymajorgrids=true,
	every axis plot/.append style={thick}
]
\addplot[color=cyan]
coordinates {(0, 676.8393589743589)
(10, 676.806282051282)
(20, 676.8460256410257)
(30, 676.7066666666666)
(40, 676.0402564102565)
(50, 674.605641025641)
(60, 671.2811538461539)
(70, 661.4074358974358)
(80, 639.4719230769231)
(90, 600.0228205128205)
(100, 547.7225641025641)};
\end{axis}

\begin{axis}[
	at={($(hybrid_alpha.east)+(2.3cm,0)$)},anchor=west,
	width=7.5cm,
	height=4cm,
	title={Performance of hybrid method vs $\tau$},
	xlabel=Planning horizon for dynamic schedule ($\tau$),
	ylabel=Average contribution,
	xmin=1, xmax=6,
	xmajorgrids=true,
	ymajorgrids=true,
	every axis plot/.append style={thick}
]
\addplot[color=cyan]
coordinates {(1, 685.068076923077)
(2, 683.4012820512821)
(3, 680.7535897435898)
(4, 662.6794871794872)
(5, 656.8829487179487)
(6, 654.390641025641)};
\end{axis}
\end{tikzpicture}
\caption{Influence of $\alpha$ and $\tau$ on performance of the hybrid method} \label{fig:hybrid_full}
\end{figure}

Figure \ref{fig:hybrid_full} shows the average contribution per time period over one year from applying this method for different values of $\alpha$ and $\tau$ in  $50$ trials. When $\alpha$ is varied, $\tau$ is fixed to $3$ time periods. When $\tau$ is varied, $\alpha$ is fixed to $50\%$. Ideally, we choose $\alpha$ and $\tau$ as high as possible without significant performance loss, as this means that more appointments can be booked further in advance. 

\colr{With $\alpha = 0\%$ being fully dynamic and $\alpha = 100\%$ fully static, we find that there is almost no performance loss when SMK statically allocates $1\%$ to $60\%$ of timeslots. For $\tau$, the steepest performance decrease happens when we go  from $\tau = 3$ to $\tau = 4$. Following this, when using the hybrid method, we recommend \colr{SMK to fix} $60\%$ of the OD appointments using the static allocation method, and allocating the remaining $40\%$ of timeslots $3$ time periods ($6$ weeks) in advance using the rolling-horizon LP. This way SMK can achieve the performance of the rolling-horizon LP for only 3 periods in advance, while planning patients up to 6 periods in advance. }

\subsubsection{Performance on access times and unused capacity}\label{sec:results_access_times}
While it is clear from the results that the solution methods improve on the static allocation method used by SMK in terms of the contribution function, we compare how well the methods adhere to the \colr{access time targets and use of capacity as the aims of this paper are to identify and exploit the effects of the allocation of timeslots on access times and use of resource capacity (see Section \ref{sec:introduction})}. We use the static allocation method as the SMK baseline measurement, and compare it to the hybrid method with $\alpha = 60\%$, $\tau = 3$, and the rolling-horizon LP with planning horizon $p=6$ (twelve weeks).

Table \ref{tab:access_times_results} shows the percentage of patients meeting their \colr{access time targets} and average access time of late patients (in time periods) for these methods, where the average is taken over $100$ trials in the simulation.

\begin{table}[h!]
\caption{Access time statistics achieved by hybrid method ($\tau=3, \alpha=60\%$), rolling-horizon LP ($p=6$), and static allocation.\\}
\centering
\label{tab:access_times_results}
\begin{tabular}{|c|ccc|ccc|}
\hline
\multirow{2}{*}{\textbf{Queue}} & \multicolumn{3}{c|}{\textbf{Appointments within access time target (\%)}} & \multicolumn{3}{c|}{\textbf{Average access time (time periods)}} \\ \cline{2-7}  & \hspace*{2.5mm} \textbf{Hybrid} \hspace*{2.5mm}     & \hspace*{2.5mm} \textbf{LP} \hspace*{2.5mm}    & \textbf{Static}     & \hspace*{2.5mm} \textbf{Hybrid} \hspace*{2.5mm}      & \hspace*{2.5mm} \textbf{LP} \hspace*{2.5mm}      & \textbf{Static}      \\ \hline
$FA_2$  & $26.87$   & $26.02$    & $3.27$   & $3.83$                & $3.83$            & $4.80$               \\ \hline
$FU_3$                          & $99.97$            & $93.01$         & $99.77$             & $2.78$                & $2.84$            & $1.19$               \\ \hline
$FU_6$                          & $100$              & $100$           & $99.87$             & $5.20$                & $5.11$            & $0.67$               \\ \hline
$FU_{12}$                       & $100$              & $100$           & $99.96$             & $7.62$                & $4.34$            & $0.23$               \\ \hline
$OR_1$                          & $96.76$            & $96.09$         & $96.52$             & $0.97$                & $0.97$            & $0.99$               \\ \hline
$OR_2$                          & $97.66$              & $97.17$         & $97.05$             & $1.85$                & $1.86$            & $1.89$               \\ \hline
$OR_4$                          & $97.15$              & $96.64$         & $97.46$             & $3.47$                & $3.51$            & $3.52$               \\ \hline
$OR_6$                          & $98.60$              & $98.16$         & $97.97$             & $1.86$                & $1.94$            & $2.04$               \\ \hline
$DA_3$                          & $100$                & $82.24$         & $95.34$             & $2.46$                & $2.53$            & $0.84$               \\ \hline
\end{tabular}
\end{table}

For most appointment types, the hybrid method marginally increases the percentage of appointments that are realized within the access time target, compared to the rolling-horizon LP and static allocation. The exception is first appointments, where an improvement from $3.27\%$ to $26.87\%$ is realized. There is also a $5$ \colr{percentage point} increase in the number of on-time discharge appointments from static allocation to the hybrid method.

It is interesting to see that, while the average access times of late patients for follow-up and discharge appointments are much shorter when using the static allocation, the hybrid method achieves a higher percentage of appointments within the access time target limit, as does the rolling-horizon LP for FU's. This is because in the contribution function used in the LP, no cost is assigned for making a patient wait for another time period when they are still within their access time target. The LP also has more flexibility in prioritizing late FA or DA patients over FU patients when necessary, and is better at handling unexpected outlier patients. This presents a trade-off to hospital management: ensure that on average, patients move through the treatment process quickly, or ensure that as many patients as possible are treated within their target? It is not immediately obvious that these goals are in conflict, yet these results show that this is the case. Sometimes, less urgent patients must wait for a relatively longer time in order to ensure that more urgent patients can be treated on time.

Another goal was to make maximal use of the available resource capacity. In the same experiment as above, we also measured how much OD and OR capacity was left unused by both methods during one year. The hybrid method left $0.76\%$ of OD capacity and $0.31\%$ of OR capacity unused on average. The rolling-horizon LP wasted $1.15\%$ of OD capacity and $0.35\%$ of OR capacity. The static allocation wasted $2.70\%$ of OD capacity, but only $0.18\%$ of OR capacity. However, during this experiment (as with all other experiments in this paper), the integer constraint on actions in the ILP was relaxed to a non-negativity constraint. Otherwise, the experiment would have taken around twenty hours to run, as opposed to the fifteen minutes it took with the relaxation. This results in some wasted capacity, as the actions selected by the LP must be rounded down to the nearest integer to ensure feasibility. We expect less wasted capacity, and even better access times, when the integer constraint is applied.

We should point out that these results are very dependent on the cost and reward parameters, in the case of the LP and hybrid method, and on the proportions used in the static allocation method. For example, increasing the proportion of first appointments in the static allocation \colr{would have resulted in better access times while having a} negative impact on FU and DA access times. As for costs and rewards, assigning a higher cost for making FA patients wait will likely lead to shorter external access times but longer internal access times. Of course, it is up to the hospital to decide how to prioritize each appointment type.

%% file: conclusion.tex
\section{Conclusions and recommendations} \label{sec:conclusion_recommendations}

In this paper, we allocate timeslots to consultation types in sessions at the outpatient department and operating room of one surgeon for a case study at the orthopaedic department of the Sint Maartenskliniek in Nijmegen. Our goal is to reduce the access time of late patients for their treatment, positively impacting their health and satisfaction. Furthermore, we want to make maximal use of the available resource capacity (i.e., OD and OR time), in order to treat as many patients as possible, and keep idle time to a minimum.

The timeslot allocation problem can be classified as a sequential decision problem under uncertainty, and was modelled as a Markov decision process. The objective is to minimize access times of late patients and maximize the number of patients treated. We also consider that appointment schedules must be generated well in advance, and develop a method to predict the future state for which planning decisions must be made.

\colr{To solve the timeslot allocation problem, we developed several solution methods, namely LSPI, a rolling-horizon LP, and several decision rules. }When comparing these solution methods, the rolling-horizon LP resulted in the highest contribution gained, with the \colr{Highest Contribution} decision rule coming in second place. Although LSPI performed best \colr{for the case study instance with 400 initial patients}, it performed worse than the other two methods when the number of patients increase, and performed worse than the \colr{currently used} static allocation method \colr{at SMK} when planning more than $10$ weeks ahead. We suspect that the increased size and complexity of the problem no longer permits a good solution in the form of a linear parametrisation.

We exploited the performance increase caused by planning less far ahead by proposing a hybrid method, where $\alpha\%$ of OD appointments are fixed in advance by a static allocation, and the remaining appointments are planned $\tau$ time periods ahead using the rolling-horizon LP. Experiments showed that $\alpha = 60\%$ and $\tau = 3$ ($6$ weeks) resulted in a high average contribution while ensuring that as many appointments as possible can be booked further in advance. This hybrid method was compared to the static allocation and rolling-horizon LP with a planning horizon of $12$ weeks to determine the actual impact on reduction of access times and unused resource capacity. The configuration of cost and reward parameters used resulted in a $23$ \colr{percentage point} increase in the number of FA appointments within the access time target, and a $5$ \colr{percentage point} increase for DA appointments, compared to the static allocation method. For other appointment types, \colr{we find} small increases in the \colr{percentage} of on-time appointments. Adjusting the cost and reward parameters changes which patients are prioritised by the LP and will therefore lead to different results. Finally, we should note that better results can be expected from the hybrid method and rolling-horizon LP when the integer constraint is applied.

We recommend using the hybrid method with either the rolling-horizon ILP or the Highest Contribution decision rule in practice, as these performed much better than LSPI on the real-life problem, are easier to implement and are very fast. Note that we use the term rolling-horizon \textit{ILP} here, as, for practical use, the integer constraint on actions should be applied in the program to ensure optimality and less wasted resource capacity. The parameters $\alpha=60\%$ and $\tau=3$ can be used with the hybrid method, but should be confirmed using data of other surgeons.

\colr{The values of the costs and rewards in the contribution function should be carefully chosen by the hospital as the policy provided by the solution methods is highly dependent on it. This could be done by repeatedly simulating the solution using different cost and reward settings as we have done, and calculating performance metrics from the simulation, until desired metrics are achieved. An iterative method to update the costs and parameters is proposed by \citet{Hulshof2013} and could also be employed. Alternatively, the costs and rewards could be derived from the actual financial costs associated with wasting resource capacity and increasing access times. We did not focus on this aspect of the process as our main goal is to compare the proposed solution methods in terms of the contribution function.}

The choice between the rolling-horizon ILP and decision rule depends on the software available to the hospital. The state prediction algorithm, decision rule, and simulation can be programmed in a standard programming language.  
However, to solve the ILP, the hospital must have access to a solver such as Gurobi. 
If this is not the case, we recommend using the decision rule, as the performance decrease is not significant, and solvers cost money. Another option is to relax the integer constraint on actions and use an open-source solver, as these are capable of solving linear programs.

Only a simplified version of our proposed method \colr{was} implemented at SMK, because full implementation would require hours of manual input and output per week to the \colr{hospital information} system. This is often a bottleneck when implementing advanced methods in practice, since most hospital systems cannot embed these methods or communicate with other systems. Because of the manual work, 80 – 100\% of the timeslots at SMK are fixed using static allocation 12 weeks in advance, which reduces the flexibility of allocating the remaining timeslots dynamically 6 weeks in advance. This simplification still led to a reduction of the waiting list for follow-up appointments with an access time target of 6 weeks while keeping the other waiting lists stable.

Although many adaptations to LSPI were attempted to improve performance and convergence, many more can be identified. However, even if we were able to improve the performance of LSPI, there are still three major drawbacks to the algorithm: complexity of implementation, running time, and explainability. First, the algorithm requires much more data collection in order to program a simulation, and is tricky to implement. Convergence is very sensitive to settings of parameters and problem size. Second, on realistic problem sizes, the algorithm takes hours to converge, making experimentation with parameters practically impossible. For practical use, the algorithm would have to be run every few months, to ensure model parameters are up-to-date. The high running time also implies that we cannot hope to expand the model to optimize for multiple surgeons simultaneously. \colr{The third and maybe} largest drawback for practical implementation is that the algorithm is a black box: we cannot reasonably explain why one action is chosen over the other. In healthcare applications, this means that we cannot explain why we accept one patient, but reject the other. This is definitely not the case for the decision rule (treat the most urgent patient), and for the ILP we can argue that the actions chosen are the optimal solution to an explainable problem.

One of the simplifying assumptions placed on the model was that a surgeon is not affected by activities of other surgeons. This allowed us to create a single-surgeon model. In real-life, this is not the case. For example, at SMK, new patients are not assigned to a single surgeon, but to a subspecialty consisting of multiple surgeons. Therefore, if one surgeon has a long \colr{FU} waiting list, another surgeon can compensate by accepting more new patients. The ILP formulated in this paper can be expanded to model a subspecialty, or even the entire orthopaedic department. This would make the program larger, and integer constraints on variables would be necessary in order to allocate whole patients to surgeons, significantly increasing solving time. However, this could be easier for the capacity planning department as the model would output \colr{timeslot allocations} for all surgeons simultaneously, rather than having to repeatedly solve the ILP for every surgeon separately. Furthermore, being able to take interaction effects into account could significantly improve performance.

The fact that the Highest Contribution decision rule overtook the rolling-horizon LP in terms of average contribution when there was a high number of initial patients in the system, suggested that the rolling-horizon LP could be improved by decreasing $\gamma$ when the number of patients increases. Further research could be done on developing an adaptive-$\gamma$ scheme to keep the LP outperforming the decision rule. The need to decrease the discount factor will likely not only rely on the number of patients, but also on the available resource capacity, \colr{as systems with more capacity can handle more patients on time. We expect that this} adjustment to the ILP will not significantly affect computing time, \colr{while improving} performance.

%% file: abbreviations.tex
\section{Abbreviations}\label{sec:abb}

\begin{itemize}
  \item ADP: Approximate Dynamic Programming
  \item DA: Discharge Appointment
  \item EVI: Exact Value Iteration
  \item FA: First Appointment
  \item FU: Follow-Up
  \item ILP: Integer Linear Program
  \item LP: Linear Program
  \item LSPI: Least-squares policy iteration
  \item LSTD: Least-squares temporal differencing
  \item MILP: Mixed Integer Linear Program
  \item MDP: Markov Decision Process  
  \item MSE: Mean Squared Error
  \item OD: Outpatient Department
  \item OR: Operating Room
  \item SMK: Sint Maartenskliniek
\end{itemize}